\documentclass[10pt]{amsart}

\usepackage{amssymb}
\usepackage{amsmath}
\usepackage{graphicx}

\newtheorem{stheorem}{Theorem}[section]
\newtheorem{slemma}[stheorem]{Lemma}
\newtheorem{sproposition}[stheorem]{Proposition}
\newtheorem{scorollary}[stheorem]{Corollary}
\newtheorem{_claim}{Claim}

\newtheorem*{maintheorem}{Main Theorem}

\newtheorem{remk}[stheorem]{Remark}
\newenvironment{sremark}{\begin{remk}\rm}{\end{remk}}

\numberwithin{equation}{section}
\numberwithin{table}{section}
\numberwithin{figure}{section}

\newenvironment{sproof}{\begin{proof}\rm}{\end{proof}}
\newenvironment{mainproof}%
{\begin{proof}[Proof of Main Theorem]\rm}{\end{proof}}

\newcommand{\C}{\mathord{\mathbb C}}
\newcommand{\F}{\mathord{\mathbb F}}
\renewcommand{\P}{\mathord{\mathbb  P}}
\newcommand{\Q}{\mathord{\mathbb  Q}}
\newcommand{\R}{\mathord{\mathbb  R}}
\newcommand{\Z}{\mathord{\mathbb Z}}

\newcommand{\FF}{\mathord{\mathcal F}}

\newcommand{\NN}{\mathord{\mathcal N}}
\newcommand{\OO}{\mathord{\mathcal O}}
\newcommand{\UU}{\mathord{\mathcal U}}
\newcommand{\XX}{\mathord{\mathcal X}}
\newcommand{\YY}{\mathord{\mathcal Y}}
\newcommand{\ZZ}{\mathord{\mathcal Z}}

\renewcommand{\aa}{\mathord{\bf a}}
\newcommand{\bb}{\mathord{\bf b}}
\newcommand{\cc}{\mathord{\bf c}}

\newcommand{\inv}{\sp{-1}}
\newcommand{\sprime}{\sp{\prime}}
\newcommand{\spprime}{\sp{\prime\prime}}
\newcommand{\rep}[1]{\langle #1 \rangle}

\newcommand{\transpose}[1]{#1\sp{T}}

\newcommand{\id}{\mathop{\rm id}}

\renewcommand{\Im}{\mathop{\rm Im}\nolimits}
\newcommand{\Ker}{\mathop{\rm Ker}\nolimits}
\newcommand{\Coker}{\mathop{\rm Coker}\nolimits}
\newcommand{\Hom}{\mathop{\rm Hom}\nolimits}

\newcommand{\Mat}{\mathop{\rm Mat}\nolimits}
\newcommand{\Card}{\mathop{\rm Card}\nolimits}
\newcommand{\PGL}{\mathop{PGL}\nolimits}

\newcommand{\Imag}{\mathop{\frak{Im}}\nolimits}
\newcommand{\Real}{\mathop{\frak{Re}}\nolimits}

\newcommand{\Spec}{\mathop{\rm Spec}}
\newcommand{\Proj}{\mathop{\rm Proj}}

\renewcommand{\mod}{\,\mathop{\rm mod}\,}

\newcommand{\binomial}[2]
{\Bigl(
\begin{array}{c}
#1 \\ #2\\
\end{array}
\Bigl) }

\newcommand{\noproof}{\hfill $\qed$}

\newcommand{\surj}{\mathbin{\to \hskip -8pt \to}}

\newcommand{\maprightsp}[1]{\smash{\mathop{\longrightarrow}\limits\sp{#1}}}

\newcommand{\maprightsb}[1]{\smash{\mathop{\longrightarrow}\limits\sb{#1}}}

\newcommand{\mapleftsp}[1]{\smash{\mathop{\longleftarrow}\limits\sp{#1}}}

\newcommand{\mapleftsb}[1]{\smash{\mathop{\longleftarrow}\limits\sb{#1}}}

\newcommand{\vstrechmapdown}{\llap{\phantom{\bigg\downarrow}}}

\newcommand{\mapdown}{
\hbox{$\bigm\downarrow$}
\vstrechmapdown
}

\newcommand{\mapdownsurj}{
\hbox{$\bigm\downarrow$}
\llap{\hbox{\raise 2pt\hbox{$\bigm\downarrow$}}}%
\vstrechmapdown
}

\newcommand{\mapdownleft}[1]{
\llap{$\vcenter{\hbox{$\scriptstyle#1$}}$}
\mapdown
}

\newcommand{\mapdownright}[1]{
\mapdown
\rlap{$\vcenter{\hbox{$\scriptstyle#1$}}$ }
}

\newcommand{\mapdownsurjleft}[1]{
\llap{$\vcenter{\hbox{\raise 6pt\hbox{$\scriptstyle#1$}}}$}
\mapdownsurj
}

\newcommand{\mapdownsurjright}[1]{
\mapdownsurj
\rlap{$\vcenter{\hbox{\raise 6pt\hbox{$\scriptstyle#1$}}}$}
}

\newcommand{\homeo}
{\setbox1=\hbox{$\to$}%
\setbox2=\hbox to \wd1 {\hfill$\scriptstyle\sim$\hskip 1.2pt\hfill}%
\,\raise 3.7pt\box2\llap{\box1}\,%
}
\newcommand{\isom}{\homeo}

\newcommand{\map}[3]{{#1} \, : \, {#2} \, \to \, {#3}}
\newcommand{\homeomap}[3]{{#1} \, : \, {#2} \,\homeo \, {#3}}
\newcommand{\mapsurj}[3]{{#1} \, : \, {#2} \, \surj \, {#3}}

\newcommand{\itcolon}{\hbox{\hskip .8pt\rm :}}

\def\Beta{B\sb{\eta}}
\def\ve{\varepsilon}
\def\spsh{\sp{\sharp}}
\def\spfl{\sp{\flat}}
\relax{
\def\bdrsh{\bdr\spsh}
\def\bdrfl{\bdr\spfl}
}
\def\bdrsh{\bdr\sp{B}}
\def\bdrfl{\bdr\sp{B}}
\def\bdr{\partial}
\def\spci{\sp{ci}}
\def\Ha{H\sb{n,\aa}}
\def\Hb{H\sb{n,\bb}}
\def\Hc{H\sb{n,\cc}}
\def\Sa{S\sb{n, \aa}}
\relax{
\def\Ha{\Hilb (n, \aa)}
\def\Hb{\Hilb (n,\bb)}
\def\Hc{\Hilb (n,\cc)}
\def\Sa{S(n, \aa)}
}
\relax{
\def\Ha{H\sb{\aa}}
\def\Hb{H\sb{\bb}}
\def\Hc{H\sb{\cc}}
\def\Sa{S\sb{\aa}}
}
\def\Ma{(M\sb{\aa})\sb 0}
\def\Na{(N\sb{\aa})\sb 0}
\def\Maci{\Ma \spci}
\def\Mb{(M\sb{\bb})\sb 0}
\def\Mbci{\Mb \spci}
\def\IMa{(I\sb o M\sb{\aa})\sb 0}
\def\IMaci{\IMa \spci}
\def\IMb{(I\sb o M\sb{\bb})\sb 0}
\def\IMbci{\IMb \spci}
\def\HM{ \Hom (M\sb{\bb}, M\sb{\aa})\sb 0}
\def\HMN{ \Hom (M\sb{\bb}, N\sb{\aa})\sb 0}
\def\Fba{F\sb{\bb, \aa}}
\def\tlZZ{\widetilde\ZZ}
\def\IN{(I\sb o N\sb{\aa})\sb 0}
\def\INci{\IN\spci}
\def\spdi{\sp{\natural}}
\def\hg{\rep {h(g)}}

\def\oghg{(o, \rep  g, \hg)}
\def\pghg{(p, \rep  g, \hg)}
\def\Xhg{X\sb{\hg}}
\def\sbgh{\sb{(g, h)}}
\def\sbaa{\sb{\aa}}
\def\sbbb{\sb{\bb}}
\def\sbcc{\sb{\cc}}
\def\Maa{M\sbaa}
\def\Mbb{M\sbbb}
\def\Mcc{M\sbcc}
\def\Naa{N\sbaa}
\def\Io{I\sb o}

\def\quand{\quad\textrm{and}\quad}
\def\sbD{\sb{\Delta}}
\def\sbo{\sb o}
\def\spp{\sp{(p)}}

\begin{document}

\setlength{\textwidth}{13.5cm}

\title[Vanishing cycles, GHC and Gr\"obner bases]{Vanishing cycles, the generalized
Hodge Conjecture and Gr\"obner bases}

\author{Ichiro Shimada}
\address{
Division of Mathematics,
Graduate School of Science,
Hokkaido University,
Sapporo 060-0810, JAPAN
}
\email{shimada@math.sci.hokudai.ac.jp}

\begin{abstract}
Let $X$ be a general complete intersection of a given multi-degree
in a complex projective space.
Suppose that the anti-canonical line bundle of $X$ is ample.
Using the cylinder homomorphism
associated with the family of complete intersections contained
in $X$,
we prove that the vanishing cycles in 
the middle homology group of $X$
are represented by topological cycles
whose support is 
contained in a proper Zariski closed subset $T\subset X$
of certain codimension.
In some cases,
we can find such a Zariski closed subset $T$
with  codimension equal to the upper bound
obtained from the Hodge structure of
the middle cohomology group
of $X$ by means of Gr\"obner bases. Hence 
a consequence of the generalized Hodge conjecture
is verified in these cases.
\end{abstract}



\maketitle

\section{Introduction}\label{sec:introduction}
There are only few non-trivial examples that can be used 
as supporting evidence 
for the generalized Hodge conjecture
formulated by Grothendieck  \cite{Grothendieck69}.
In this paper, we deal with complete intersections  of
small multi-degrees in a complex projective
space,
and prove,
in some cases,
a consequence of the generalized Hodge conjecture
for  these complete intersections
by means of cylinder homomorphisms.
\par
%
We work over the complex number field $\C$.
Let $X$ be a general complete intersection
of multi-degree $\aa=(a\sb 1, \dots, a\sb r)$ in $\P\sp n$
with $\min (\aa)\ge 2$.
Suppose that  $X$ is Fano, that is,
the total degree 
 $\sum\sb{i=1}\sp r a\sb i$ of $X$ is  less than or equal to $n$.
We put 
$$
m:=\dim X=n-r
\quand
k:= \Bigl[\frac{1 }{\max (\aa)} \Bigl( n-\sum\sb{i=1}\sp r a\sb i \Bigr) \Bigr]+1, 
$$
where $[\phantom{a}]$ denotes the integer part.
It is known 
that the Hodge structure of the middle cohomology group $H\sp m (X,\Q)$  of 
$X$ satisfies the following (\cite[Expos\'e XI, Corollaire 2.8]{SGA7II}):
\begin{equation}\label{eq:HS}
H\sp{\nu,m-\nu} (X)=0 
\quad \Longleftrightarrow \quad
0\le \nu < k
\;\textrm{ or }\;
0\le m-\nu < k. 
\end{equation}
If the generalized Hodge conjecture is true, then 
there should  exist a Zariski closed subset
$T$ of $X$ with codimension $k$
such that the inclusion $T\hookrightarrow X$
induces a surjective homomorphism
$H\sb m (T, \Q)\surj H\sb m (X, \Q)$.
\par
We will try to verify this consequence of
the generalized Hodge conjecture
by means of  {\em cylinder homomorphisms}.
Let $\bb=(b\sb 1, \dots , b\sb s)$
be another sequence of
 integers
satisfying 
$\min (\bb)\ge 1$ and $r< s<n$.
We denote by $F\sb{\bb} (X)$ the  scheme 
parameterizing  all complete intersections of multi-degree $\bb$ in $\P\sp n$
that are contained in $X$,
and by $Z\sb{\bb} (X)\subset X\times  F\sb{\bb} (X)$ the universal family
with
\begin{equation}\label{diag:fam1}
\begin{array}{ccc}
Z\sb{\bb} (X) &\maprightsp{\alpha\sb X} &X \\
\mapdownleft{\pi\sb X} && \\
F\sb{\bb}(X) && \\
\end{array}
\end{equation}
being the diagram of the projections.
We put $l:=n-s$.
Suppose that $F\sb{\bb} (X)$ is non-empty,
and  that  
$m>2l$ holds.
Since $\pi\sb X$ is proper and flat of relative dimension $l$, 
we have a homomorphism
$$
H\sb{m-2l} (F\sb{\bb} (X) , \Z) \,\to\, 
 H\sb{m} (Z\sb{\bb} (X) , \Z)
$$
that maps a homology class
$[\tau]\in H\sb{m-2l} (F\sb{\bb} (X) , \Z)$
represented by a topological $(m-2l)$-cycle $\tau$
in $F\sb{\bb} (X)$
to the homology class
$[\pi\sb X \inv (\tau)]\in H\sb m (Z\sb{\bb} (X), \Z)$
represented by the topological $m$-cycle 
$\pi\sb X \inv (\tau)$
in $Z\sb{\bb} (X)$.
We define a homomorphism
$$
\map{\psi\sb{\bb} (X)}{ H\sb{m-2l} (F\sb {\bb} (X), \Z) }{ H\sb m (X, \Z)}
$$
by $\psi\sb{\bb} (X) ([\tau]):=\alpha\sb{X*} ([\pi\sb X \inv (\tau)])$,
and 
call  $\psi\sb{\bb} (X)$
{\em the cylinder homomorphism associated with  the family
$\pi\sb X : Z\sb{\bb} (X)\to F\sb{\bb}(X)$}.
\par
\smallskip
It was remarked in \cite{Steenbrink85}
that there exists a Zariski closed subset $T$ of $X$
with codimension $\ge l$ such that
the image of the homomorphism
$H\sb m (T, \Q)\to H\sb m (X, \Q)$
induced from the inclusion $T\hookrightarrow X$
contains $\Im \psi\sb{\bb} (X)\otimes \Q$.
(See also Corollary~\ref{cor:toGHC}
 of this paper.)
Therefore,
in view of the generalized Hodge conjecture,
it is an interesting problem to find
a sequence $\bb$ with $l$ as large as possible
(hopefully $l=k$)
such that
the cylinder homomorphism 
$\psi\sb{\bb} (X)$
has a ``big" image.
\par
\smallskip
Our Main Theorem,
which will be stated in \S\ref{sec:statementMT},
gives us a sufficient condition on
$(n, \aa,\bb)$ for the image of 
$\psi\sb{\bb} (X)$
to contain the {\em module of vanishing cycles}
$$
V\sb m (X, \Z) :=\Ker(H\sb m (X, \Z)\to H\sb m (\P\sp n, \Z)).
$$
This sufficient condition
can be checked by means of Gr\"obner bases.
Combining Main Theorem
with a theorem of Debarre and Manivel~\cite[Th\'eor\`eme~2.1]{DebarreManivel98}
about the variety of linear subspaces contained in a general complete intersection,
we  also give a simple numerical 
condition on $(n, \aa, \bb)$
that is sufficient for $\Im \psi\sb{\bb} (X)\supseteq V\sb m (X, \Z)$
to hold (Theorem~\ref{thm:numlin}).
In many cases, 
our method yields
$\bb$ with $l$ larger than any previously known results,
and sometimes
we can verify the consequence of the generalized Hodge conjecture.
See \S\ref{sec:GHC} for the examples.
\par
After the work of Clemens and Griffiths
\cite{ClemensGriffiths72}
on the family of lines in a cubic threefold,
many authors have studied the cylinder homomorphisms of  type
$\psi\sb{\bb} (X)$,
and proved that the image contains the vanishing cycles
(%
\cite{BeauvilleDonagi85},
\cite{Clemens83},
\cite{Collino86},
\cite{Letizia84},
\cite{Lewis85},
\cite{Lewis88},
\cite{PiccoBotta89},
\cite{Shimada90a},
\cite{Shimada90b},
\cite{Shimada91},
\cite{Terasoma90},
\cite{Welters81}%
).
Our method
provides us with   a unified proof
and a generalization 
of  these  results.
\par
%
%
This paper is organized as follows.
In \S\ref{sec:statementMT},
we state Main Theorem.
In \S\ref{sec:vccyl},
we study a connection between
vanishing cycles and cylinder homomorphisms
in  general setting.
Theorem~\ref{thm:vccyl} in this section is essentially same as the result  of \cite{Shimada95}.
However we present a complete and improved  proof for  readers' convenience.
In \S\ref{sec:univ},
we construct the universal family of
the families $Z\sb{\bb} (X)\to F\sb{\bb} (X)$
over the  scheme parameterizing  all complete intersections
of multi-degree $\aa$ in $\P\sp n$,
which is a Zariski open subset of a Hilbert scheme,
and studies its properties.
Combining the results in \S\ref{sec:vccyl} and \S\ref{sec:univ},
we prove Main Theorem in
\S\ref{sec:proofMT}.
In \S\ref{sec:GBmethod},
we explain a method for checking  the conditions on $(n, \aa, \bb)$  required by Main Theorem
by means of  Gr\"obner bases.
In \S\ref{sec:DMmethod},
an application of the theorem of Debarre and Manivel is presented.
Examples  are investigated 
in  relation to the generalized Hodge conjecture in \S\ref{sec:GHC}.
\par
\smallskip
\noindent
{\bf Conventions.}
(1)
We work over  $\C$.
A point of a scheme  means a $\C$-valued point
unless otherwise stated.
(2)
For an analytic space $X$ or a scheme $X$ over $\C$,
let $T\sb p X$ denote the Zariski tangent space to $X$
at a point $p$ of $X$. 
(3)
The multi-degree of a complete intersection is always denoted in the {\em non-decreasing} order.
\section{Statement of Main Theorem}\label{sec:statementMT}
We fix an integer $n\ge 4$.
Let 
$\aa=(a\sb 1, \dots, a\sb r)$ and $\bb=(b\sb 1, \dots, b\sb s)$ be sequences of integers
satisfying
\begin{equation}\label{eq:numaabb}
2\le a\sb 1\le\dots \le a\sb r, \quad
1\le b\sb 1\le \dots \le b\sb s \quand 
r< s<n.
\end{equation}
We put
$$
m:=n-r \quand l:=n-s.
$$
We denote by $\Ha$
the  scheme parameterizing all complete
intersections of multi-degree $\aa$ in $\P\sp n$.
For a point $t$ of $\Ha$,
we denote by $X\sb t$ the corresponding complete intersection.
Let $\Sa$ denote the Zariski closed subset of  $\Ha$
parameterizing all singular complete intersections.
It is well-known that
$\Sa$ is an irreducible hypersurface of $\Ha$,
and that,
if $u$ is a general point of $\Sa$,
then $X\sb u$ has only one singular point $p$.
For $t\in \Ha$, 
we denote by
$F\sb{\bb} (X\sb t)$ the  scheme parameterizing
all complete intersections of multi-degree $\bb$ in $\P\sp n$
that are contained in $X\sb t$ as  subschemes.
If $m>2l$ and $F\sb{\bb} (X\sb t)\ne \emptyset$, then  we have the cylinder homomorphism
$$
\map{\psi\sb{\bb} (X\sb t)}{ H\sb{m-2l} (F\sb{\bb} (X\sb t), \Z)}{ H\sb m (X\sb t, \Z)}.
$$
%
%
We put
$
t\sb a :=\Card \{ \: i \mid a\sb i=a\sb r\:\}
$ and $
t\sb b :=\Card \{\: j \mid b\sb j=a\sb r\:\}.
$
%
%
%
%
%
%
\begin{maintheorem}
Suppose that the following inequalities are satisfied\itcolon
\begin{eqnarray}
& a\sb i \ge b\sb i\quad (i=1, \dots, r),
\quad
a\sb r \ge b\sb s, \quand  \label{eq:num1}\\
&m-2l \ge t\sb b - t\sb a,
\quad
m>2l.
\label{eq:num2}
\end{eqnarray}
Suppose also that,
for a general point $u$ of $\Sa$,
there exists a complete intersection
of multi-degree $\bb$ in $\P\sp n$
that is contained in $X\sb u$,
passing through the  unique singular point $p$ of $X\sb u$,
and smooth at $p$.
Then,
for a general point $t$ of $\Ha$,
the scheme $F\sb{\bb} (X\sb t)$ is non-empty, and 
the image of the cylinder homomorphism
$\psi\sb{\bb} (X\sb t)$
contains the module of vanishing cycles $V\sb m (X\sb t, \Z)$.
\end{maintheorem}
\begin{sremark}
In Proposition~\ref{prop:equivconds},
we will give
several  conditions
equivalent to the second condition of Main Theorem.
One of them can be tested easily
by means of Gr\"obner bases,
as will be explained in \S\ref{sec:GBmethod}.
\end{sremark}
\section{Vanishing cycles and a cylinder homomorphism}\label{sec:vccyl}
In this section,
we work in the category of
complex analytic spaces and holomorphic maps.
We study in  general setting the problem when the image of a cylinder homomorphism
contains a given vanishing cycle.
For the detail of the classical theory of  vanishing cycles,
we refer to \cite{Lamotke81}.
\par
%
%
Let $\varphi : Y\to\Delta$ be a proper surjective
holomorphic map from a smooth irreducible
complex analytic space
of dimension $m+1\ge 2$
to the open unit disk $\Delta\subset \C$.
For a point $a\in \Delta$,
we denote by $Y\sb a$ the fiber $\varphi\inv (a)$.
Suppose that $\varphi$ has only one critical point $p$,
that $p$ is on the central fiber $Y\sb 0$,
and that the Hessian
$$
\map{H}{T\sb p Y\times T\sb p Y}{\C}
$$
of $\varphi$ at $p$ is non-degenerate.
We put $\Delta\sp{\times}:=\Delta\setminus\{0\}$.
For any $\varepsilon \in \Delta\sp{\times}$,
the kernel of the homomorphism
$H\sb m (Y\sb{\varepsilon}, \Z)\to H\sb m (Y, \Z)$
induced from the inclusion $Y\sb{\varepsilon}\hookrightarrow Y$
is generated by the
{\em vanishing cycle}
$[\Sigma\sb{\varepsilon}]\in H\sb m (Y\sb{\varepsilon}, \Z)$
associated to the non-degenerate critical point $p$ of $\varphi$.
\par
Let $\varrho :F\to \Delta$ be a surjective  holomorphic map 
from a smooth irreducible complex analytic space
$F$ of dimension $k$  to the unit disk,
and let $W$ be a reduced closed analytic
subspace of $Y \times\sb{\Delta} F$
such that
the  projection
$\varpi : W\to F$
is flat of relative dimension $l>0$.
Since $\varphi$ is proper, so is  $\varpi$.
Let $\gamma : W\to Y$ be the projection onto the first factor.
We obtain the following commutative diagram:
\begin{equation}
\label{diagram:W}
\begin{array}{ccc}
W & \maprightsp{\gamma} & Y\\
\mapdownleft{\varpi}& & \mapdownright{\varphi} \\
F &\maprightsb{\varrho}&\Delta .\\
\end{array}
\end{equation}
For  $u\in F$,
the fiber  $\varpi\inv (u)$ can be regarded as 
a closed  $l$-dimensional analytic subspace of $Y\sb{\varrho (u)}$ by $\gamma$.
For $a\in \Delta$, we put
$F\sb a :=\varrho\inv (a)$ and 
$W\sb a:=\varpi\inv (F\sb a)$.
Then we obtain a family of $l$-dimensional closed analytic subspaces of $Y\sb a$:
\begin{equation}
\label{diagram:Wa}
\begin{array}{ccc}
W\sb a & \maprightsp{} & Y\sb a\\
\mapdownleft{\varpi\sb a}& &  \\
F\sb a.  &&\\
\end{array}
\end{equation}
Since  the restriction $\varpi\sb{a}: W\sb a  \to F\sb a$ of $\varpi$ 
to $W\sb a$
is proper and flat of relative dimension $l$,
we have the cylinder homomorphism
$$
\map{\psi\sb a}{H\sb{m-2l} (F\sb a, \Z)}{ H\sb m (Y\sb a, \Z)}
$$
associated with  the family \eqref{diagram:Wa} for any $a\in \Delta$.
\begin{stheorem}\label{thm:vccyl}
We assume  $m>2l>0$.
\par
{\rm (1)}
Suppose that there exists a point $q$ of $W\sb 0$ 
such that $\gamma (q)$ is the critical point $p$ of $\varphi$,
that $\varpi$ is smooth at $q$, and that $\gamma$ is an immersion at $q$.
Then $k=\dim F$ is less than or equal to $m-2l+1$.
\par
{\rm (2)}
Suppose moreover that $k=m-2l+1$ holds.
Then $\varpi (q)$ is a critical point of $\varrho$,
and the Hessian of $\varrho$ at $\varpi (q)$ is non-degenerate.
Let $\varepsilon$ 
be a point of $\Delta\sp{\times}$ with $|\ve|$ small enough,
and let $[\sigma\sb{\varepsilon}]\in H\sb{m-2l} (F\sb{\varepsilon}, \Z)$
be the vanishing cycle associated
to the non-degenerate critical point $\varpi (q)$ of $\varrho$.
If the vanishing cycle $[\Sigma\sb{\varepsilon}]\in H\sb m (Y\sb{\varepsilon}, \Z)$
is not a torsion element,
then  $\psi\sb{\varepsilon} ([\sigma\sb{\varepsilon}])$ 
is equal to  $[\Sigma\sb{\varepsilon}]$ up to sign.
\end{stheorem}
\begin{sproof}
(1)  Let $U\sb{W, q}$ be a small open connected neighborhood of $q$ in $W$.
We can   assume that $\varpi$ is smooth at every point
of $U\sb{W,q}$,
and  that $\gamma$ embeds $U\sb{W, q}$ into $Y$.
We put
$$
o:=\varpi (q) \quand Z:=\varpi\inv (o).
$$
Then  $\gamma (U\sb{W, q}\cap Z)$ and $\gamma(U\sb{W, q})$
are smooth locally closed analytic subsets of $Y$ passing through $p$.
Let $T\sb 1$ and $T\sb 2$ be the Zariski tangent spaces
to $\gamma (U\sb{W, q}\cap Z)$ and $\gamma(U\sb{W, q})$ at $p$,
respectively.
We have $T\sb 1\subseteq T\sb 2 \subseteq T\sb p Y$
and
$\dim T\sb 1 =l$, $\dim T\sb 2 =k+l$.
We will show that $T\sb 1$ and  $T\sb 2$ are orthogonal with respect to the Hessian
$H$ of $\varphi$ at $p$.
Let $v$ be an arbitrary vector of $T\sb 1$.
Since the structure
$\varpi | U\sb{W, q} : U\sb{W, q} \to F$
of the smooth fibration on $U\sb{W, q}$
is carried over to $\gamma (U\sb{W, q})$,
there exists a holomorphic vector field $\tilde v$ 
defined in a small open  neighborhood $U\sb{Y, p}$ of $p$ in $Y$ such that
$\tilde v\sb p$ is equal to $v$, and that,
if $q\sprime \in U\sb{W, q}$ satisfies $\gamma (q\sprime) \in U\sb{Y, p}$,
then $\tilde v\sb{\gamma (q\sprime)}$ is tangent to the smooth locally closed analytic subset 
$\gamma (U\sb{W, q}\cap \varpi\inv (\varpi (q\sprime)))$ of $Y$.
Since the diagram~\eqref{diagram:W} is commutative,
the function $\varphi$ is constant on $\gamma (\varpi\inv (\varpi (q\sprime)))$
for any $q\sprime \in U\sb{W, q}$,
and hence 
the holomorphic function $\tilde v (\varphi)$ is
constantly zero on $\gamma (U\sb{W, q} )\cap U\sb{Y, p}$, 
which means that the following holds for any $w\in T\sb 2$:
$$
H(w,v):=w (\tilde v (\varphi))=0.
$$
Thus  $T\sb 1$ is contained in the orthogonal complement $T\sb 2\sp{\perp}$ of $T\sb 2$
with respect to $H$.
Since $H$ is non-degenerate,
we have
$$
l=\dim T\sb 1 \le \dim T\sb p Y - \dim T\sb 2 =(m+1)-(k+l).
$$
Therefore we obtain  $k\le m+1-2 l$.
\par
(2)
From now on,
 we assume  $k=m+1-2l$.
Then we have $T\sb 1 =T\sb 2\sp{\perp}$.
Hence $H$ induces a non-degenerate
symmetric bilinear form 
$$
\map{H\sprime }{T\sb 2/ T\sb 1\times T\sb 2/T\sb 1}{\C}.
$$
Since $\varpi$ is smooth at $q$,
there is a local holomorphic section $s : U\sb{F, o} \to W$
of $\varpi$
defined in a small open neighborhood
$U\sb{F, o}$ of $o=\varpi (q)$ in $F$
such that $s(o)=q$.
We take $U\sb{F, o}$ so small that 
$s(U\sb{F, o}) \subset U\sb{W, q}$ holds.
Let $S$ be the image of $\gamma\circ s$,
which is a smooth locally closed analytic subset of $Y$
passing through  $p$, 
and let $T\sb 3$ be the Zariski tangent space to $S$ at $p$.
%
%
We have
$T\sb 2 =T\sb1\oplus T\sb 3$.
It follows from the non-degeneracy of $H\sprime$
that the restriction $H|T\sb 3: T\sb 3 \times T\sb 3 \to \C$
of $H$ to $T\sb 3$
is also non-degenerate.
Since
$\gamma\circ s$ yields an isomorphism from
$U\sb{F, o}$ to $S$,
and 
$\varrho$ coincides on $U\sb{F, o}$ with
$$
U\sb{F, o} \;\maprightsp{\gamma\circ s}\; S \;\maprightsp{\varphi|S}\; \Delta,
$$
the point $o$ is a critical point of $\varrho$.
Moreover, 
the Hessian of $\varrho$ at $o$ is equal to
$H|T\sb 3$
via the isomorphism $ (d \;(\gamma\circ s))\sb o : T\sb o F \;\isom\; T\sb 3$,
and hence is non-degenerate.
\par
We will describe the holomorphic maps in the diagram~\eqref{diagram:W} 
in terms of local coordinates.
Let $t$ be the coordinate on $\Delta$.
There exist  local analytic coordinates 
$x=(x\sb 1, \dots, x\sb k)$ 
on $F$ with the center $o$
such that $\varrho$ is given by
\begin{equation}\label{eq:rhot}
\varrho\sp * t = x\sb 1\sp 2 + \cdots + x\sb k \sp 2.
\end{equation}
Since $\varpi$ is smooth at $q$,
there exists  a local analytic coordinate system 
$(w, w\sprime)=(w\sb 1, \dots, w\sb k, w\sprime\sb 1, \dots, w\sprime\sb l)$ 
on $W$ with the center $q$
such that $\varpi$ is given by 
\begin{equation}\label{eq:pix}
\varpi\sp * x\sb i = w\sb i \qquad (i=1, \dots, k).
\end{equation}
Since $\gamma$ is an immersion  at $q$,
there exist  local analytic coordinates 
$(y, y\sprime, y\spprime)=
(y\sb 1, \dots, y\sb k, y\sprime\sb 1, \dots, y\sprime\sb l,y\spprime\sb 1,\dots,
y\spprime\sb l)$ 
on $Y$ with the center $p$
such that $\gamma$ is given by
\begin{equation}\label{eq:gammay}
\begin{cases}
\gamma\sp * y\sb i  = w\sb i & (i=1, \dots, k),\\
\gamma\sp * y\sprime\sb j  = w\sprime\sb j & (j=1, \dots, l),\\
\gamma\sp * y\spprime\sb j  = 0  &(j=1, \dots, l).\\
\end{cases}
\end{equation}
(Note that $\dim Y$ is equal to $m+1=k+2l$.)
Then the locally closed analytic subset $\gamma (U\sb{W, q})$ of
$Y$
is defined by $y\spprime\sb1=\dots=y\spprime\sb l=0$
locally around   $p$.
From the commutativity of the diagram~\eqref{diagram:W},
it follows that
 $\varphi\sp * t$
and 
$y\sb 1\sp 2 + \cdots +y\sb k\sp 2$
coincide on  $\gamma (U\sb{W, q})$.
Therefore,
in a small neighborhood of $p$,
the function $\varphi\sp * t$ is written as follows:
$$
y\sb 1\sp 2 + \cdots +y\sb k\sp 2+ a\sb 1 y\spprime\sb 1 + \cdots +a\sb l y\spprime\sb l,  
$$
where $a=(a\sb1, \dots, a\sb l)$ is a system of  holomorphic functions defined 
locally around $p$.
Since $p$ is a critical point of $\varphi$,
we have
$a\sb 1 (p) =\dots =a\sb l (p)=0$.
The non-degeneracy of the Hessian $H$ of $\varphi$ at $p$
implies that the $l\times l$ matrix  $\bigl({\partial a\sb i }/{\partial y\sprime\sb j} (p)\bigr)\sb{i,j=1, \dots, l}$
is non-degenerate.
Hence
$(y, a, y\spprime)$ is another local analytic coordinate system
on $Y$ with the center $p$.
We replace $y\sprime$ with $a$.
Then  we have
\begin{equation}\label{eq:phit}
\varphi\sp * t = 
y\sb 1\sp 2 + \cdots +y\sb k\sp 2+ 
y\sprime\sb 1 y\spprime\sb 1 + \cdots +y\sprime\sb l
y\spprime\sb l.
\end{equation}
We can make coordinate transformation
on $w\sprime$
according to the coordinate transformation on $y\sprime$
so that \eqref{eq:gammay} remains valid.
We put
\begin{equation}\label{eq:ytoz}
\begin{cases}
z\sb i  := y\sb i & (i=1, \dots, k),\\
z\sb{k+j}  := (y\sprime \sb j + y\spprime\sb j)/2 & (j=1, \dots, l),\\
z\sb{k+l+j}  := \sqrt{-1}\;(y\sprime\sb j-y\spprime\sb j)/2  &(j=1, \dots, l).\\
\end{cases}
\end{equation}
Then we have
\begin{equation}\label{eq:phitz}
\varphi\sp * t = 
z\sb 1\sp 2 + \cdots +z\sb {m+1}\sp 2.
\end{equation}
\par
Let $\eta$ be a sufficiently small positive real number,
and let $\Beta$ be the closed ball in $Y$ defined by
$$
|z\sb 1|\sp 2 + \cdots +|z\sb{m+1}|\sp 2 \le \eta.
$$
Let $\ve$ be a positive real number
that is small enough compared with $\eta$.
Let $s$ be a real number satisfying $0<s\le \ve$.
The closed subset
$$
Y\sb s \,\cap\, \Beta  \,=\,\{\,(z\sb 1, \dots, z\sb{m+1})\,\mid\,
|z\sb 1|\sp 2 + \cdots +|z\sb{m+1}|\sp 2 \le \eta,\,
 z\sb 1\sp 2 + \cdots + z\sb{m+1}\sp 2=s
\,\}
$$ 
of $Y\sb s=\varphi\inv (s)$
is homeomorphic to the total space
$$
E: =\{\;(u, v) \in \R\sp{m+1}\times \R\sp{m+1}\;\mid\;
\| u\|=1,\; \| v\|\le 1, \;u\perp v
\;\}
$$
of the unit disk tangent bundle $\tau: E\to S\sp m$
of the $m$-dimensional sphere 
$S\sp m : =\{u \in \R\sp{m+1}\mid\| u\|=1\}$,
where the projection $\tau$ is given by $\tau(u, v)=u$.
We identify $S\sp m$ with the zero section of $\tau: E\to S\sp m$.
The homeomorphism
$h\sb s : Y\sb s \cap \Beta \homeo E$
is written explicitly as follows:
\begin{equation}\label{eq:h}
u=\frac{\Real (z)}{\|\Real (z)\|},
\qquad
v=\sqrt{\frac{2}{\eta-s}}\Imag (z).
\end{equation}
Its inverse $h\sb s\inv : E \homeo Y\sb s \cap \Beta $ is given by the following:
\begin{equation}\label{eq:hinv}
z=\sqrt{s+\Bigl( \frac{\eta-s}{2}\Bigr) \| v\|\sp 2}\cdot  u +
\sqrt{-\Bigl(\frac{\eta-s}{2}\Bigr)}\cdot v.
\end{equation}
The sphere $S\sp m \subset E$ is mapped by $h\sb s\inv$ to
the closed submanifold
$$
\Sigma\sb s :=
\left\{ \;(z\sb 1, \dots, z\sb{m+1})\in Y\;\biggm\vert\; 
\vcenter{\hbox{\vbox{
\hbox{$ z\sb 1\sp 2 + \cdots + z\sb{m+1}\sp 2=s$,}
\vskip 3pt
\hbox{$\Imag (z\sb i)=0 \;\; (i=1, \dots, m+1)$}
}
}
}
\;\right\}
$$
of $Y\sb s$.
With an orientation,
this topological  $m$-cycle
$\Sigma\sb s$ represents the vanishing cycle
$[\Sigma\sb s ]\in H\sb m (Y\sb s, \Z)$,
which generates the kernel of the homomorphism 
$H\sb m(Y\sb s, \Z)\to H\sb m (Y, \Z)$
induced from  $Y\sb s \hookrightarrow Y$.
\par
%
%
For $s\in (0,\ve]$,
let $\sigma\sb s$
denote the $(m-2l)$-dimensional sphere
contained in $ F\sb s= \varrho\inv (s)$ defined by
$$
\sigma\sb s :=
\left\{ \;(x\sb 1, \dots, x\sb{k})\in F\;\biggm\vert\; 
\vcenter{\hbox{\vbox{
\hbox{$x\sb 1\sp 2 + \cdots + x\sb{k}\sp 2=s$,}
\vskip 3pt
\hbox{$\Imag (x\sb i)=0 \;\; (i=1, \dots, k)$}
}
}
}
\;\right\}.
$$
With an orientation,
this topological $(m-2l)$-cycle $\sigma\sb s$
represents
the vanishing cycle
$[\sigma\sb s]\in H\sb{m-2l} (F\sb s, \Z)$
associated to the non-degenerate critical point
$o$ of $\varrho$.
Since $\varpi$ is proper and flat of relative dimension $l$,
the inverse image $\varpi\inv (\sigma\sb{s})$
of the oriented sphere
$\sigma\sb s$ 
can be considered as a topological $m$-cycle in $W\sb s=\varpi\inv (F\sb s)$.
The image $\psi\sb{\ve} ([\sigma\sb\ve])$
of $[\sigma\sb\ve]\in H\sb{m-2l} (F\sb{\ve} , \Z)$ by the cylinder homomorphism
$\psi\sb\ve : H\sb{m-2l} (F\sb\ve, \Z)\to H\sb m (Y\sb\ve, \Z)$
is represented by the topological  $m$-cycle 
$$
\map{\gamma | \varpi\inv (\sigma\sb\ve)}{\varpi\inv (\sigma\sb\ve) }{Y\sb\ve}.
$$
Since  the sphere $\sigma\sb\ve$ bounds
an $(m-2l+1)$-dimensional  closed ball in $F$,
the topological  $m$-cycle 
$\gamma | \varpi\inv (\sigma\sb\ve)$
is a boundary
of a topological   $(m+1)$-chain in $Y$;
that is,
$\psi\sb{\ve} ([\sigma\sb\ve])$ belongs to
 the kernel of $H\sb m (Y\sb{\ve}, \Z)\to H\sb m (Y, \Z)$.
Hence
there exists an integer $c$
such that the following holds in $H\sb m (Y\sb\ve, \Z)$:
\begin{equation}\label{eq:c}
\psi\sb{\ve} ([\sigma\sb\ve])=c\, [\Sigma\sb\ve].
\end{equation}
We will show that,
if $[\Sigma\sb\ve]$ is not a torsion element in 
$H\sb m (Y\sb\ve, \Z)$, then  $c$ is $\pm 1$.
\par
%
%
%
%
We put
$$
Y\sb{[0, \ve]} :=\varphi\inv ([0, \ve])=\bigcup\sb{s\in[0, \ve]} Y\sb s.
$$
For any closed subset $A$ of $Y\sb{[0, \ve]}$,
we set
$$
A\spsh := A\setminus (A\cap \Beta\sp{\circ}),
\quad
A\spfl :=  A\cap\Beta
\quand
\bdrfl A := A\cap \bdr \Beta,
$$
%
%
%
%
where $\Beta\sp{\circ}$ is the interior of the closed ball $\Beta$,
and $\bdr \Beta$ is the boundary of  $\Beta$.
The sharp $\sharp$ means ``outside the ball",
and the flat $\flat$ means ``inside the ball".
The explicit descriptions~\eqref{eq:h}~and~\eqref{eq:hinv}
of the homeomorphism $h\sb s : Y\sb s\spfl \homeo  E$
for $s\in (0,\ve]$ show that the  restriction 
$h\sb s | \bdrfl Y\sb s : \bdrfl Y\sb s \homeo \bdr E$
of $h\sb s$ to $\bdrfl Y\sb s$ can be extended 
to a homeomorphism from
$$
\bdrsh Y\sb 0=\{\;(z\sb 1, \dots, z\sb{m+1})\;\mid\;
|z\sb 1|\sp 2 + \cdots +|z\sb{m+1}|\sp 2 = \eta,\;
 z\sb 1\sp 2 + \cdots + z\sb{m+1}\sp 2=0
\;\}
$$
to
$\bdr E=\{(u, v) \in E\mid
\| v\|= 1\}$
smoothly.
We denote these homeomorphisms by 
$$
\homeomap{\bdrsh h\sb s}{ \bdrsh Y\sb s}{\bdr E} \quad (s\in [0, \ve]).
$$
The homeomorphism
$\bdrsh h\sb 0 : \bdrsh Y\sb 0\homeo \bdr E$
is given by the following:
%
%
$$
u=\sqrt{{2}/{\eta}} \;\Real (z), \quad
v=\sqrt{{2}/{\eta}} \;\Imag (z),
\quand
z=\sqrt{{\eta}/{2}} \: (u+\sqrt{-1} v ).
$$
Putting these homeomorphisms $\bdrsh h\sb s$ $(s\in [0,\ve])$ together,
we obtain a trivialization
$$
\homeomap{\bdrsh h}{\bdrsh  Y\sb{[0, \ve]} }{\bdr E \times [0, \ve]}
$$
of the restriction
$\varphi| \bdrsh Y\sb{[0, \ve]}  : 
\bdrsh Y\sb{[0, \ve]} \to [0,\ve]$
of $\varphi$ to $\bdrsh Y\sb{[0, \ve]} $ over $[0, \ve]$.  
Let
$$
\homeomap{\bdrsh f}{\bdrsh Y\sb{[0, \ve]}}{\bdrsh Y\sb\ve \times [0, \ve]}
$$
be the trivialization of $\varphi| \bdrsh Y\sb{[0, \ve]}$
obtained by composing $\bdrsh h$ and $(\bdrsh h\sb\ve \times \id)\inv$.
Since
the only critical point $p$ of $\varphi$ is not contained in $Y\sb{[0, \ve]}\spsh$,
we can show by
Ehresmann's fibration theorem for the manifolds
with boundaries 
that
the trivialization $\bdrsh f$ extends to
a trivialization
\begin{equation}\label{eq:trivialization}
\homeomap{(f\spsh, \bdrsh f)}{(Y\sb{[0, \ve]}\spsh, \bdrsh Y\sb{[0, \ve]})}%
{(Y\sb\ve\spsh, \bdrsh Y\sb\ve)\times [0,\ve]}
\end{equation}
of $\varphi |Y\sb{[0,\ve]}\spsh: Y\sb{[0,\ve]}\spsh\to [0, \ve]$
in such a way that
the restriction 
%
%
of $(f\spsh, \bdrsh f) $ to the fiber over $\ve$
is the identity map.
For $s\in [0,\ve]$,
let
$$
\homeomap{(f\sb s\spsh, \bdrsh f\sb s)}{(Y\sb{s}\spsh, \bdrsh Y\sb{s})}{(Y\sb\ve\spsh, \bdrsh Y\sb\ve)}
$$
denote the restriction of $(f\spsh, \bdrsh f) $ to the fiber over $s$.
\par
%
%
We put
$$
C\sb s :=\gamma (\varpi\inv (\sigma\sb s)) \subset Y\sb s.
$$
When $s$ approaches $0$,
this closed subset $C\sb s$ degenerates into $C\sb 0:=\gamma (\varpi\inv (o))$,
which is an $l$-dimensional closed analytic subset of  $Y\sb 0$.
We decompose $\varpi\inv(\sigma\sb s)$
into the union of $\varpi\inv(\sigma\sb s)\sp{(\sharp)}$ and $\varpi\inv(\sigma\sb s)\sp{(\flat)}$,
where
\begin{eqnarray*}
\varpi\inv(\sigma\sb s)\sp{(\sharp)}&:=&
\varpi\inv(\sigma\sb s)
\;\setminus\; ( \gamma\inv (\Beta\sp{\circ}) \cap \varpi\inv(\sigma\sb s))
\quand  \\
\varpi\inv(\sigma\sb s)\sp{(\flat)}&:=&\gamma\inv (\Beta)  \cap \varpi\inv (\sigma\sb s).
\end{eqnarray*}
Since $\eta$ and $\ve$ are small enough,
and $W$ is a subspace of $Y\times F$,
we have
\begin{equation}\label{eq:Vsprime}
\varpi\inv(\sigma\sb s)  \sp{(\flat)}= W\cap (\Beta \times \sigma\sb s) \;\;\subset\;\; U\sb{W, q}
\end{equation}
for all $s\in [0, \ve]$,
where $U\sb{W, q}$ is the open neighborhood of $q$ in $W$
that was introduced at the beginning of the proof.
Recalling that $\gamma$ embeds $U\sb{W, q}$ into $Y$,
we see that 
the map $\gamma$ yields a homeomorphism from 
$\varpi\inv(\sigma\sb\ve) \sp{(\flat)}$
to $C\sb\ve\spfl$.
By definition, 
$\gamma$ maps 
$\varpi\inv(\sigma\sb\ve)\sp{(\sharp)}$
to $C\sb\ve\spsh$.
We then define 
a closed subset $\widetilde C\sb\ve$ of $Y\sb\ve$ by
\begin{equation}\label{eq:decompwtC}
\widetilde C\sb\ve
				\;:=\;
				C\sb\ve\spfl
							\;\cup\;
				\Bigl(\bigcup\sb{s\in [0, \ve]} \bdrsh f\sb s (\bdrsh C\sb s)  \Bigr)
						 \;\cup\; 
				f\sb 0\spsh (C\sb 0\spsh).
\end{equation}
Note that we have
$$
\bdrsh\widetilde C\sb\ve =\bigcup\nolimits\sb{s\in [0, \ve]} \bdrsh f\sb s (\bdrsh C\sb s)
\quand
\widetilde C\sb\ve\spfl = C\sb\ve\spfl \cup \bdrfl\widetilde C\sb\ve,
\quad
\widetilde C\sb\ve\spsh = \bdrsh\widetilde C\sb\ve \cup f\sb 0\spsh (C\sb 0\spsh).
$$
%
%
%
%
%
%
\par
%
%
Using the trivialization $(f\spsh, \bdrsh f)$,
we can ``squeeze" the topological  $m$-cycle 
$\gamma | \varpi\inv(\sigma\sb\ve) : \varpi\inv(\sigma\sb\ve) \to Y\sb\ve$
outside the ball
so that the image is contained in $\widetilde C\sb{\ve}$.
More precisely, 
we can construct a homotopy 
from 
$\gamma | \varpi\inv(\sigma\sb\ve) : \varpi\inv(\sigma\sb\ve) \to Y\sb\ve$
to a continuous map
$\beta : \varpi\inv(\sigma\sb\ve) \to Y\sb\ve$
with the following properties:
\par
\smallskip
\noindent
($\beta$-1)
The  image  $\beta(\varpi\inv(\sigma\sb\ve))$ of $\beta$
coincides with $\widetilde C\sb\ve$.
\par
\noindent
($\beta$-2)
The homotopy is 
stationary on 
$\varpi\inv(\sigma\sb\ve) \sp{(\flat)}$.
In particular, $\beta$ yields a homeomorphism from 
$\varpi\inv(\sigma\sb\ve)  \sp{(\flat)}$
to the first piece  $C\sb\ve\spfl$ of the decomposition~\eqref{eq:decompwtC}.
\par
\noindent
($\beta$-3)
The image  $\beta (\varpi\inv(\sigma\sb\ve)  \sp{(\sharp)})$
of $\varpi\inv(\sigma\sb\ve)  \sp{(\sharp)}$ by $\beta$ is contained in 
$\widetilde C\sb\ve\spsh$.
\par
\smallskip
\noindent
(See Figure~\ref{fig:homotopy}.)
Since  $\psi\sb\ve ([\sigma\sb\ve])$
is represented by  $\gamma | \varpi\inv(\sigma\sb\ve)$,
it is also represented by  the topological $m$-cycle $\beta$.
\begin{figure}[t]
  \centering\includegraphics{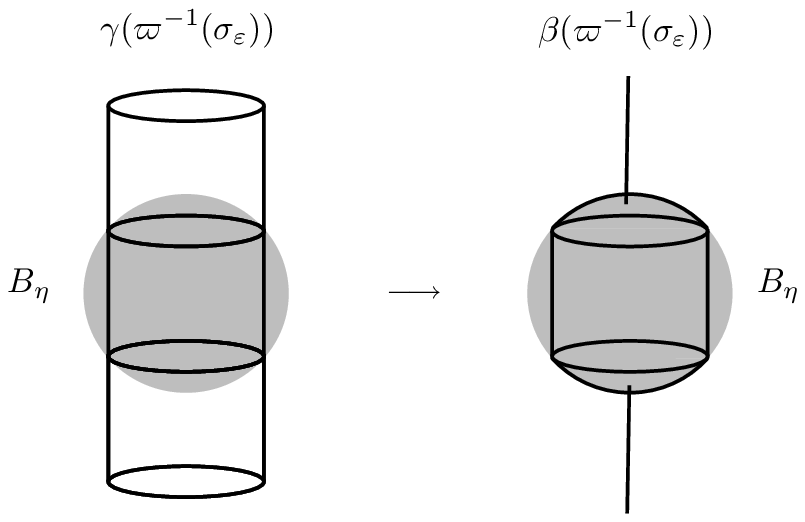}
		\caption[homotopy]{Homotopy from $\gamma | \varpi\inv(\sigma\sb\ve)$ to $\beta$}
  \label{fig:homotopy}
\end{figure}
\par
%
%
From~\eqref{eq:Vsprime} and 
\eqref{eq:rhot}, 
\eqref{eq:pix}, 
\eqref{eq:gammay},
\eqref{eq:ytoz},
we see that 
$C\sb{s}\spfl$ $(s\in [0, \ve])$  is given in terms of
the local coordinate system $z$ by the following:
$$
\begin{cases}
|z\sb 1|\sp 2 + \cdots +|z\sb{m+1}|\sp 2 \le \eta, & \\
z\sb 1\sp 2 + \cdots +z\sb{k}\sp 2=s, &\\
\Imag (z\sb i)=0 \quad (i=1, \dots, k), &\\
y\spprime\sb j=z\sb{k+j}+\sqrt{-1} z\sb{k+l+j}=0 \quad (j=1, \dots, l). &
\end{cases}
$$
For $s\in [0, \ve]$,
let 
$G\sb s$ be the closed subset of $E$ defined by the following equations:
$$
\begin{cases}
(\, 2s +(\eta-s) \|v\|\sp 2\,)(\, u\sb 1\sp 2+\cdots + u\sb k\sp 2\,)=2s, & \\
v\sb 1=\dots =v\sb k=0, &\\
v\sb{k+j}=-g\sb s (\|v\|) \cdot u\sb{k+l+j} \quad (j=1, \dots, l),& \\
v\sb{k+l+j}=g\sb s (\|v\|)\cdot  u\sb{k+j}\quad (j=1, \dots, l), &
\end{cases}
$$
where 
$$
g\sb s (\|v\|):=\sqrt{\frac{2s}{\eta-s}+\|v\|\sp 2}.
$$
Then, for
$s\in (0, \ve]$,
the homeomorphism $h\sb s : Y\sb s\spfl \homeo E$ maps 
$C\sb s\spfl$ to  $G\sb s$.
In particular,
the first piece 
$C\sb\ve\spfl$ of the decomposition~\eqref{eq:decompwtC}
of $\widetilde C\sb\ve$ is
mapped homeomorphically to $G\sb\ve$ by $h\sb\ve$.
It is easy to check that, 
for any $s\in [0, \ve]$ (including $s=0$),
the homeomorphism
$\bdrsh h\sb s : \bdrsh Y\sb s \homeo \bdr E$
maps 
$\bdrfl C\sb s$
to
$G\sb s\cap \bdr E$.
We put
$$
T\sb\ve :=\{\; u\in S\sp m \mid u\sb1 \sp 2+\cdots + u\sb k\sp 2 < 2\ve/(\eta+\ve) \;\}, 
$$
and let  $T\sb\ve\sp{-}$  be the closure of $T\sb\ve$.
We can easily check
that
the projection $\tau: E\to S\sp m$ induces a homeomorphism
from $G\sb\ve$ to $S\sp m\setminus T\sb\ve$,
and that,
for any $s\in [0, \ve]$,
$G\sb s \cap \bdr E$ is contained in $\tau\inv (T\sb\ve\sp{-})\cap \bdr E$.
In particular,
the second piece
$\bdrfl \widetilde C\sb\ve$
of the decomposition~\eqref{eq:decompwtC}
is mapped by $\bdrsh h\sb\ve$ into 
$\tau\inv (T\sb\ve\sp{-})\cap \bdr E$.
\par
%
Let $a$ be a point of $S\sp m\setminus T\sb\ve\sp{-}$.
Then the closed subset
$h\sb\ve \inv (\tau\inv (a))$ of $Y\sb\ve\spfl$
intersects
the first piece
$C\sb\ve\spfl$
of the decomposition~\eqref{eq:decompwtC}
at only one point,
which is in the interior
of $C\sb\ve\spfl$,
and the intersection is transverse.
Moreover, $h\sb\ve \inv (\tau\inv (a))$
is disjoint
from the second piece
$\bdrfl \widetilde C\sb\ve$
of the decomposition~\eqref{eq:decompwtC}.
The third piece
$f\sb 0\spsh (C\sb 0\spsh)$
is a topological  $2l$-cycle
in $(Y\sb\ve\spsh, \bdrsh Y\sb\ve)$,
because $C\sb 0$ is a
topological $2l$-cycle in $Y\sb 0$.
\par
\smallskip
If $[\Sigma\sb\ve]\in H\sb m (Y\sb\ve, \Z)$ is zero,
then 
$\psi\sb\ve ([\sigma\sb\ve])=0$
by~\eqref{eq:c} and hence there is nothing to prove.
%
%
%
%
Suppose that $[\Sigma\sb\ve]$ is not zero and 
not a torsion element. Then 
there exists a homology class $[\Theta] \in H\sb m (Y\sb\ve, \Z)$
such that the intersection number $[\Sigma\sb\ve]\cdot [\Theta]$
of $[\Sigma\sb\ve]$ and $[\Theta]$
in  $Y\sb\ve$
is not zero.
In order to show that the integer $c$ in \eqref{eq:c} is $\pm 1$,
it is enough to prove the following:
\begin{equation}\label{eq:intsame}
\psi\sb\ve ([\sigma\sb\ve])\cdot [\Theta]=\pm[\Sigma\sb\ve]\cdot [\Theta].
\end{equation}
Multiplying $[\Theta]$ by a positive integer if necessary,
we can assume that
$[\Theta]$ is represented by
a compact oriented  $m$-dimensional differentiable submanifold
$\Theta$ of $Y\sb\ve$
(\cite{Thom54}).
By the elementary transversality theorem
(see, for example, \cite{GolubitskyGuillemin73}),
we can move $\Theta$  in $Y\sb\ve$ in such a way that the following hold:
\par
\smallskip
\noindent
($\Theta$-1)
			The closed subset $h\sb\ve (\Theta\spfl)$ of $E$ is a union of finite number of
			fibers of $\tau: E\to S\sp m$ over points in $S\sp m \setminus T\sb\ve\sp{-}$.
\par
\noindent
($\Theta$-2)
			The topological $m$-cycle $\Theta\spsh$ of $(Y\sb\ve\spsh, \bdrsh Y\sb\ve)$
			is disjoint from 
			the topological  $2l$-cycle $f\sb 0\spsh (C\sb 0\spsh)$.
   Here we use the assumption $m>2l$.
\par
\smallskip
\noindent
From ($\Theta$-1) and ($\Theta$-2), 
the points $\Theta \cap \widetilde C\sb{\ve}$
are contained in the interior
of the first piece $C\sb{\ve}\spfl$
of the decomposition~\eqref{eq:decompwtC} of $\widetilde C\sb{\ve}$,
and the intersections are all transverse.
Moreover,
the total intersection number
of $\Theta$ and $\widetilde C\sb{\ve}$
is equal to
that of $\Theta$ and $\Sigma\sb{\ve}$
up to sign,
because both of them
are equal,  up to sign,  to
the number of fibers of $\tau$ constituting $h\sb\ve (\Theta\spfl)$
(counted with signs according to the orientation).
Combining these with  the properties ($\beta$-1)-($\beta$-3) of the topological $m$-cycle $\beta$,
we see that
$[\beta]\cdot [\Theta]=\pm [\Sigma\sb{\ve}]\cdot [\Theta]$.
We have seen that  $\psi\sb\ve ([\sigma\sb\ve])$
is represented by  $\beta$.
Thus we obtain~\eqref{eq:intsame}.
\end{sproof}
\section{The universal family}\label{sec:univ}
In this section, we will construct the universal family of the incidence
varieties of complete intersections in a complex projective space $\P\sp n$.
\par
%
%
First we fix some notation.
Let 
$$
R=\bigoplus\sb{d=0}\sp{\infty} R\sb d:=\C[x\sb 0, \dots, x\sb n]
$$
be the polynomial ring of $n+1$ variables
with coefficients in $\C$ graded by the degree $d$ of polynomials.
We set $R\sb d:=0$ for $d<0$.
Let $M$ be a graded $R$-module.
We denote by $M\sb d$ the vector space
consisting of homogeneous elements of $M$ with degree $d$. 
For an integer $k$,
let $M(k)$ be the $R$-module $M$ with grading
shifted by 
$M(k)\sb d:=M\sb{k+d}$.
For another graded $R$-module $N$,
let  $\Hom (M, N)\sb 0$ denote  the vector space
of degree-preserving homomorphisms  from $M$ to $N$.
Let $\cc=(c\sb 1, \dots, c\sb t)$
be a sequence of positive integers.
We assume $t<n$.
Let us  define the graded free $R$-module $M\sb{\cc}$
by
$$
M\sb{\cc}:=\bigoplus\sb{i=1}\sp t R(c\sb i).
$$
An element of $M\sb{\cc}$ is written  as a column vector.
Let $f=\transpose{(f\sb 1, \dots, f\sb t)}$ be an element
of $(M\sb{\cc})\sb 0=\oplus R\sb{c\sb i}$,
where $f\sb i$ is a homogeneous polynomial
of degree $c\sb i$.
We denote by $J\sb f$ the homogeneous ideal of $R$
generated by $f\sb1, \dots, f\sb t$.
There exists a Zariski open dense subset
$(M\sb{\cc})\sb 0\spci$
of the vector space $(M\sb{\cc})\sb 0$ consisting
of all $f\in (M\sb{\cc})\sb 0$ such that
the ideal $J\sb f$ defines a complete intersection of
multi-degree $\cc$ in $\P\sp n =\Proj R$.
For $f\in (M\sb{\cc})\sb 0\spci$,
let  $Y\sb{\rep f}$ denote the complete intersection
defined by $J\sb f$.
It is well-known that,
for any integer $\nu$,
the dimension of the vector space
$$
H\sp 0 (Y\sb{\rep f}, \OO(\nu))= (( R/J\sb f) (\nu))\sb 0
$$
is independent of the choice of $f\in (M\sb{\cc})\sb 0\spci$.
\par
%
%
Let $\Hc$ denote  the  scheme parameterizing all complete intersections
of  multi-degree $\cc$ in $\P\sp n$.
It is well-known that  $\Hc$ is  a smooth irreducible  quasi-projective scheme.
For an element $f\in (M\sb{\cc})\sb 0\spci$,
let $\rep f$ denote the point of $\Hc$
corresponding to the complete intersection $Y\sb{\rep f}$.
We have a  surjective morphism
$$
\mapsurj{q\sb{\cc}}{ (M\sb{\cc})\sb 0\spci}{ \Hc}
$$
that maps $f$ to $\rep f$.
Let $\YY\sb{\cc}\subset \P\sp n \times \Hc$ be the universal family of complete intersections of
multi-degree $\cc$ in $\P\sp n$
with $\phi\sb{\cc} : \YY\sb{\cc} \to \Hc$ and 
$\tau\sb{\cc} : \YY\sb{\cc} \to \P\sp n$ the projections.
\begin{sproposition}
{\rm (1)} The morphism $q\sb{\cc}$ is smooth.
{\rm (2)} The morphism $\tau\sb{\cc}$ is smooth.
In particular, $\YY\sb{\cc}$ is smooth.
\end{sproposition}
\begin{sproof}
(1)
The Zariski tangent space
to $\Hc$ at $\rep f$ is
given by
\begin{equation}\label{eq:THilb}
T\sb{\rep f} \Hc = 
H\sp 0 (Y\sb{\rep f}, \NN\sb{Y\sb{\rep f}/\P\sp n}) =
(M\sb{\cc}/J\sb f M\sb{\cc})\sb 0,
\end{equation}
where $\NN\sb{Y\sb{\rep f}/\P\sp n}$ is the normal sheaf
of $Y\sb{\rep f}$ in $\P\sp n$,
which is isomorphic to $\oplus\sb{i=1}\sp t \OO (c\sb i)$.
By~\eqref{eq:THilb} and 
$T\sb f (M\sb{\cc})\sb 0\spci \cong (M\sb{\cc})\sb 0$,
the linear map $(d q\sb{\cc})\sb f :T\sb f (M\sb{\cc})\sb 0\spci \to T\sb{\rep f} \Hc$
is identified with the quotient  homomorphism
$(M\sb{\cc})\sb 0 \surj (M\sb{\cc}/J\sb f M\sb{\cc})\sb 0$.
Hence $q\sb{\cc}$ is smooth.
\par
(2)
Let $P=(p, \rep{f})$ be a point of $\YY\sb{\cc}$,
where $p$ is a point of $Y\sb{\rep{f}}$,
and let $I\sb p$ be the homogeneous ideal
of $R$ defining the point $p$.
The kernel of $(d\tau\sb{\cc})\sb{P} : T\sb{P} \YY\sb{\cc} \to T\sb p \P\sp n$
is mapped isomorphically to a subspace
of $T\sb{\rep{f}} \Hc$ by
$(d\phi\sb{\cc})\sb{P} : T\sb{P} \YY\sb{\cc} \to T\sb{\rep{f}} \Hc$.
This subspace coincides with the subspace 
$(I\sb p \Mcc/ J\sb{f} \Mcc)\sb 0$ 
of $(\Mcc/ J\sb{f} \Mcc)\sb 0$ under the identification~\eqref{eq:THilb}.
Since 
$\dim (\Mcc / I\sb p \Mcc)\sb 0=t$,
 $\dim \Ker (d\tau\sb{\cc})\sb{P}$ is 
 equal to $\dim \Hc -t=\dim \YY\sb{\cc} -n$
for any point $P\in \YY\sb{\cc}$.
\end{sproof}
\par
%
%
Let $\aa=(a\sb 1, \dots, a\sb r)$
and $\bb=(b\sb1, \dots, b\sb s)$
be  two sequences of integers
satisfying~\eqref{eq:numaabb}.
Instead of $\YY\sb{\aa}$ and $\YY\sb{\bb}$, 
we denote by
\begin{equation}\label{diag:univsXZ}
\begin{array}{ccc}
\XX &\maprightsp{\tau} & \P\sp n\\
\mapdownleft{\phi} && \\
\Ha &&\phantom{\Ha} \\
\end{array}
\quand\quad
\begin{array}{ccc}
\ZZ &\maprightsp{\tau\sprime} & \P\sp n\\
\mapdownleft{\phi\sprime} && \\
\Hb&&\phantom{\Hb} \\
\end{array}
\end{equation}
the universal families
over $\Ha$ and $\Hb$.
%
%
For $f\in \Maci$ and $g\in \Mbci$,
we denote
by $X\sb{\rep f}$ and 
$Z\sb{\rep g}$ the complete intersections
corresponding to
$\rep f \in \Ha$ and $\rep g \in \Hb$,
respectively.
\par
An element
$h$ of $\HM$ is expressed by an $r\times s$ matrix
$(h\sb{ij})$ with
$h\sb{ij}\in R\sb{a\sb i-b\sb j}$.
When $g \in \Mb$ is fixed,
the image of the linear map
$\Hom (\Mbb, \Maa)\sb 0 \to\Ma$
given by $h\mapsto h(g)$
coincides with $(J\sb g \Maa)\sb 0$.
The  following proposition is then obvious:
\begin{sproposition}\label{prop:three}
The following  three conditions
on  the pair  $(f, g)$ of $f\in \Maci$ and $g\in \Mbci$ are equivalent\itcolon
\par
{\rm (i)}
    $X\sb{\rep f}$ contains $Z\sb{\rep g}$ as a subscheme,
\par
{\rm (ii)}
    $f$ is contained in $(J\sb g\Maa)\sb 0$, and 
\par
{\rm (iii)}
    there exists an element $h\in \HM$ such that $f=h(g)$.
\noproof
\end{sproposition}
\par
Let $\FF\sb{\bb, \aa}$ be the contravariant functor
from the category
of locally noetherian schemes
over $\C$ to the category of sets
that associates
to a locally noetherian scheme $S\to\Spec\C$
the set of pairs
$(Z\sb S, X\sb S)$,
where
$Z\sb S \subset \P\sp n\times S$ and $X\sb S \subset \P\sp n\times S$
are  families of complete intersections in $\P\sp n$ with multi-degrees $\bb$ and $\aa$,
respectively,
parameterized by $S$
such that $Z\sb S$ is a subscheme of $X\sb S$.
This functor $\FF\sb{\bb,\aa}$
is represented by a closed subscheme $\Fba$
of  $\Hb\times\Ha$.
(The scheme $\Fba$ may possibly be empty.)
We denote the projections by
$\rho\sprime : \Fba \to \Hb$ and $\rho : \Fba \to \Ha$. 
The universal family over $\Fba$ is the pair $(\tlZZ, \widetilde \XX)$
of
$\tlZZ:=\ZZ\times\sb{\Hb}\Fba$ and
$\widetilde \XX:=\XX\times\sb{\Ha}\Fba$.
We denote by
$\pi: \tlZZ \to\Fba$ and $\beta : \tlZZ \to \ZZ$
the natural projections.
We also denote by $\alpha :\tlZZ\to \XX$
the composite of the closed immersion
$\tlZZ\hookrightarrow \widetilde \XX$
and the natural projection $\widetilde \XX\to \XX$.
Thus we obtain the following commutative diagram:
\begin{equation}\label{diag:univ}
\begin{array}{ccccccccc}
\P\sp n & \mapleftsp{\tau\sprime}\hskip 0pt  & \ZZ & \mapleftsp{\beta} & \tlZZ &\maprightsp{\alpha} 
& \XX &\hskip 0pt\maprightsp{\tau} &
\P\sp n
\\
\phantom{\Hb}&&\mapdownleft{\phi\sprime} &\square &\mapdownright{\pi} &&\mapdownright{\phi}&& \\
&&\Hb &\mapleftsb{\rho\sprime} &\Fba &\maprightsb{\rho}&\Ha,&&
\phantom{\Ha} \\
\end{array}
\end{equation}
in which $\tau\circ\alpha=\tau\sprime\circ\beta$ holds.
A point of $\tlZZ$ is a triple
$$
(p, \rep g,\rep f)\in \P\sp n \times \Hb\times \Ha
$$
that satisfies  $p\in Z\sb{\rep g}\subset X\sb{\rep f}$.
The projection $\pi$ maps $(p, \rep g,\rep f)$ to
$(\rep g, \rep f)\in \Fba$,
and the morphism $\alpha$ maps $(p, \rep g,\rep f)$ to $(p, \rep f)\in \XX$.
\par
The right square 
of the diagram~\eqref{diag:univ}
is the universal family of the families~\eqref{diag:fam1} of complete intersections of multi-degree $\bb$
contained in  complete intersections of multi-degree $\aa$.
Remark that
the linear automorphism group $\PGL (n+1)$ of $\P\sp n$
acts on the diagram~\eqref{diag:univ}.
\begin{sremark}
If $(n,\aa,\bb)$ satisfies the first inequality $a\sb i\ge b\sb i$ $(i=1, \dots, r)$
of the condition~\eqref{eq:num1}
in Main Theorem,
then $\Fba$ is non-empty.
Indeed, we  choose
linear forms $\ell\sb 1, \dots, \ell\sb r, \ell\sb 1\sprime, \dots, \ell\sb s\sprime \in R\sb 1$
generally.
We define $g\in \Mbci$ by $g\sb j :={\ell\sb j\sprime}\sp{ b\sb j}$.
Since $a\sb i \ge b\sb i$,
we can define $f\in \Maci$ by
$f\sb i := {\ell\sb i\sprime}\sp{b\sb i} {\ell\sb i}\sp{ a\sb i -b\sb i}$.
Then $(\rep{g}, \rep{f})$ is a point of $\Fba$.
\end{sremark}
From now on to the end of this section, we assume that $\Fba$ is non-empty.
We define a vector space $U$ 
with a natural morphism
$\nu : U\to \Ma$ by
$$
U:=\Mb\times \HM\quand \nu (g, h):=h(g).
$$
We then put
$$
U\spci:=(\Mbci\times\HM)\cap \nu\inv (\Maci).
$$
Note that $U\spci$ is a Zariski open subset of $U$,
and hence is irreducible.
By Proposition~\ref{prop:three},
the map
$$
\sigma (g, h) :=(\rep g, \rep{h(g)})
$$
defines a surjective morphism
$\sigma : U\spci\surj \Fba$,
which makes the following diagram commutative:
$$
\begin{array}{ccc}
U\spci & \maprightsp{\nu|U\spci} &\Maci \\
\mapdownsurjleft{\sigma} && \mapdownsurjright{q\sb{\aa}} \\
\Fba &\maprightsb{\rho} &\Ha.
\end{array}
$$
In particular,
the scheme $\Fba$ is irreducible.
%
%
%
%
\begin{sproposition}\label{prop:rhosprime}
The morphism $\rho\sprime: \Fba\to\Hb$ is smooth.
\end{sproposition}
\begin{sproof}
For a non-negative integer $k$,
we set $A\sb k :=\C[ t ]/ (t\sp{k+1})$,
and for a scheme $T$ over $\C$,
we denote by $T(A\sb k)$ the set of $A\sb{k}$-valued points of $T$.
Suppose  we are given 
$\rep g\sp{[k+1]}\in\Hb (A\sb{k+1})$ and
$(\rep g\sp{[k]}, \rep f\sp{[k]})\in \Fba (A\sb k)$ 
satisfying
$\rep g \sp{[k]}=\rep g\sp{[k+1]}\mod t\sp{k+1}$.
It is enough to show that
$(\rep g\sp{[k]}, \rep f\sp{[k]})$
extends to an
element $(\rep g\sp{[k+1]}, \rep f\sp{[k+1]})$ of $\Fba (A\sb{k+1})$
over the given point 
$\rep g\sp{[k+1]}\in \Hb (A\sb{k+1})$.
Since 
$q\sb{\aa} :\Maci\surj\Ha$ and
$q\sb{\bb} : \Mbci\surj\Hb$ are smooth,
there exist 
$$
g\sp{[k+1]} \in \Mb\otimes\sb{\C} A\sb{k+1}\quand  
f\sp{[k]} \in \Ma\otimes\sb{\C} A\sb{k}
$$
 that satisfy  $q\sb{\bb} (g\sp{[k+1]})=\rep g\sp{[k+1]}$
and $q\sb{\aa} (f\sp{[k]})=\rep f\sp{[k]}$.
We put
$$
g\sp{[k]} := g\sp{[k+1]}\mod t\sp{k+1}\in \Mb\otimes\sb{\C} A\sb k,
$$
which satisfies $\rep{g\sp{[k]}}=\rep g \sp{[k]}$.
By the definition of $\Fba$,
the ideal $J\sb{g\sp{[k]}}$ of $R\otimes\sb{\C}A\sb{k}$
generated by the components of $g\sp{[k]}$ contains
the ideal $J\sb{f\sp{[k]}}$.
Hence there exists 
$h\sp{[k]}\in \HM\otimes\sb{\C} A\sb{k}$
such that $f\sp{[k]}=h\sp{[k]}(g\sp{[k]})$ holds.
Let 
$h\sp{[k+1]}$
be any element of  $\HM\otimes\sb{\C} A\sb{k+1}$ satisfying 
$h\sp{[k+1]}\mod t\sp{k+1} =h\sp{[k]}$.
We put
$$
f\sp{[k+1]} := h\sp{[k+1]} (g\sp{[k+1]} ) \in \Ma\otimes\sb{\C} A\sb{k+1}.
$$
Since being a complete intersection is an open condition
on  defining polynomials,
the ideal $J\sb{f\sp{[k+1]}}$ of $R\otimes\sb{\C}A\sb{k+1}$
defines a family of complete intersections 
of multi-degree $\aa$ over $\Spec A\sb{k+1}$.
Thus $(\rep{g\sp{[k+1]}}, \rep{f\sp{[k+1]}})$ is the hoped-for $A\sb{k+1}$-valued point of  $\Fba$.
\end{sproof}
\begin{scorollary}\label{cor:FbsandtlZZmooth}
{\rm (1)}
The scheme $\Fba$ is smooth.
{\rm (2)}
The morphism $\beta:\tlZZ \to\ZZ$ is smooth.
In particular,    $\tlZZ$ is smooth.
\noproof
\end{scorollary}
Let $(g, h)$ be a point of $U\spci$.
We have the following natural identifications
of vector spaces:
\begin{eqnarray}
\label{eq:ZNZ}
&&H\sp 0 (Z\sb{\rep g}, \NN\sb{Z\sb{\rep g}/\P\sp n}) = 
T\sb{\rep g} \Hb =(\Mbb / J\sb g \Mbb)\sb 0,  \\
\label{eq:ZNX}
&&H\sp 0 (Z\sb{\rep g}, \NN\sb{\Xhg/\P\sp n}|Z\sb{\rep g}) =
 (\Maa / J\sb g \Maa)\sb 0,  \\
\label{eq:XNX}
&&H\sp 0 (\Xhg, \NN\sb{\Xhg/\P\sp n}) =
T\sb{\hg} \Ha =(\Maa / J\sb {h(g)} \Maa)\sb 0.\phantom{aaaaa}  
\end{eqnarray}
The restriction map 
$\NN\sb{\Xhg/\P\sp n}\to \NN\sb{\Xhg/\P\sp n}|Z\sb{\rep g} $ 
of coherent sheaves
induces, via~\eqref{eq:XNX},
a linear map
$$
\map{\zeta\sprime}{T\sb{\rep{h(g)}} \Ha}{H\sp 0 (Z\sb{\rep g}, \NN\sb{\Xhg/\P\sp n}|Z\sb{\rep g})}.
$$
Under the identifications~\eqref{eq:XNX} and~\eqref{eq:ZNX},
the linear map $\zeta\sprime$ is identified
with the natural quotient homomorphism
$$
(\Maa / J\sb {h(g)} \Maa)\sb 0 \;\surj\; (\Maa / J\sb g \Maa)\sb 0.
$$
In particular, $\zeta\sprime$  is surjective.
On the other hand,
since $Z\sb{\rep g}$ is a subscheme of $\Xhg$, there is a natural 
homomorphism 
$$
 \NN\sb{Z\sb{\rep g}/\P\sp n} \;\to\; \NN\sb{\Xhg/\P\sp n}|Z\sb{\rep g}
$$
of coherent sheaves over $Z\sb{\rep g}$,
which induces, via~\eqref{eq:ZNZ}, a linear map
$$
\map{\zeta}{T\sb{\rep g} \Hb}{H\sp 0 (Z\sb{\rep g}, \NN\sb{\Xhg/\P\sp n}|Z\sb{\rep g})}.
$$
Under the identifications~\eqref{eq:ZNZ} and~\eqref{eq:ZNX},
the linear map $\zeta$ is identified
with the   homomorphism
$$
\map{\rep{h}\sb g}{ (\Mbb / J\sb g \Mbb)\sb 0 }{(\Maa / J\sb g \Maa)\sb 0}
$$ 
induced from  $h : \Mbb \to \Maa$.
\begin{sproposition}\label{prop:diagT}
Let $(g, h)$ be a point of $U\spci$, and 
$P$  the point $\sigma (g, h)=(\rep{g}, \rep{h(g)})$ of $\Fba$.
Then we have the following diagram
of fiber product\itcolon
\begin{equation}\label{diag:T}
\begin{array}{ccc}
T\sb{P} \Fba & \maprightsp{(d\rho)\sb{P}} &T\sb{\hg} \Ha\\
\mapdownsurjleft{(d\rho\sprime)\sb{P}} & \square & \mapdownsurjright{\zeta\sprime} \\
T\sb{\rep g} \Hb &\maprightsb{\zeta} & H\sp 0 (Z\sb{\rep g}, \NN\sb{\Xhg/\P\sp n}|Z\sb{\rep g}).\\
\end{array}
\end{equation}
\end{sproposition}
\begin{sproof}
By the identifications~\eqref{eq:ZNZ} and~\eqref{eq:XNX},
any vectors of $T\sb{\rep{g}} \Hb$
and 
$T\sb{\rep{h(g)}} \Ha$
are given as  elements
$$
\bar g\sprime :=g\sprime \mod (J\sb{g} \Mbb)\sb 0
\quand
\bar f\sprime :=f\sprime \mod (J\sb{h(g)} \Maa)\sb 0
$$
of $(\Mbb/J\sb{g} \Mbb)\sb 0$
and $(\Maa/J\sb{h(g)} \Maa)\sb 0$
by some $g\sprime \in \Mb$
and $f\sprime \in \Ma$, respectively.
Let $\ve$ be a dual number: $\ve\sp 2 =0$.
The vectors
$\bar g\sprime$
and
$\bar f\sprime$
correspond
to the  infinitesimal displacements
$$
{Z}\sb{\rep {g+\ve g\sprime}} \;\to\; \Spec \C [\ve] 
\quand
{X}\sb{\rep {h(g)+\ve f\sprime}} \;\to\; \Spec \C [\ve] 
$$
of $Z\sb{\rep g}$ and  $X\sb{\rep {h(g)}}$ 
defined by the homogeneous ideals
$J\sb g + \ve J\sb{g\sprime}$ and
$J\sb {h(g)} + \ve J\sb{f\sprime}$ of $R \otimes\sb{\C} \C [\ve]$,
respectively.
Then the vector $(\bar g\sprime, \bar f\sprime)$,
regarded as 
a tangent vector to $\Hb\times\Ha$ at $(\rep{g}, \rep{h(g)})$,
is tangent to  $\Fba$
if and only if
${Z}\sb{\rep {g+\ve g\sprime}} $ is
contained in ${X}\sb{\rep {h(g)+\ve f\sprime}} $
as a subscheme;
that is, there exist
elements $h\sb1,  h\sb 2 \in \HM$
such that
$h\sb1+ \ve h\sb 2 \in \HM\otimes\sb{\C} \C[\ve]$
satisfies
the following:
\begin{equation}\label{eq:h1h2}
(h\sb1+ \ve h\sb 2 )(g+\ve g\sprime ) = h(g) +\ve f\sprime.
\end{equation}
Suppose that  $h\sb1+ \ve h\sb 2 $
satisfies~\eqref{eq:h1h2}.
Because $h\sb 1 (g)=h(g)$,
each row vector
of the matrix $h\sb 1-h$ is contained in the syzygy
of the regular sequence $(g\sb 1, \dots, g\sb s)$,
and hence every component of $h-h\sb 1$
is contained in $J\sb g$.
Therefore
the two linear maps
$\rep h \sb g$
and
$\rep{h\sb 1} \sb g$
from $(\Mbb/J\sb g \Mbb)\sb 0$ to $(\Maa/J\sb g \Maa)\sb 0$
are the same.
The equality  $h\sb 1 (g\sprime)+h\sb 2 (g)=f\sprime$
then tells us that
$f\sprime \mod (J\sb g \Maa)\sb 0$
is equal to $\rep h\sb g (\bar g\sprime)$,
because $h\sb 2 (g) \in (J\sb g \Maa)\sb 0$.
Hence $(\bar g\sprime, \bar f\sprime)$
is contained in the fiber product of $\zeta$ and $\zeta\sprime$.
Conversely,
if $(\bar g\sprime, \bar f\sprime)$
is contained in the fiber product of $\zeta$ and $\zeta\sprime$,
then it is easy to find $h\sb 2 \in \HM$
satisfying
$(h+ \ve h\sb 2 )(g+\ve g\sprime ) = h(g) +\ve f\sprime$.
\end{sproof}
Since  $\Fba$ is reduced by Corollary~\ref{cor:FbsandtlZZmooth} (1), we obtain the following:
\begin{scorollary}\label{cor:dimFbaanddimCokerdrho}
Let $(g, h)$ be an arbitrary point of $U\spci$.
\par
{\rm (1)}
The dimension of $\Fba$ is equal to
\begin{equation}\label{eq:dimFba}
\begin{split}
&\dim (\Maa/J\sb{h(g)} \Maa)\sb 0 +
\dim(\Mbb/J\sb{g} \Mbb)\sb 0-
\dim(\Maa/J\sb{g}\Maa)\sb 0
=\\
&\qquad\qquad \dim \Ha +
\dim \Hb-
\dim(\Maa/J\sb{g}\Maa)\sb 0. 
\end{split}
\end{equation}
\par
{\rm (2)}
Let $P$ be the point $\sigma (g,h)$ of $\Fba$. Then
the dimension of the cokernel of 
$(d\rho)\sb{P}: T\sb {P}\Fba \to T\sb{\rep {h(g)}} \Ha$
is equal to
$$
\hbox to \textwidth{
\hfill
$\dim \Coker\zeta 
=\dim \Coker \rep h \sb g
=\dim (\Maa/(J\sb g \Maa + h(\Mbb)))\sb 0.$
\hfill
\qed
}
$$
\end{scorollary}
\begin{sproposition}\label{prop:dimKerdalpha}
Let $(g, h)$ be a point of $U\spci$,
and let 
$p$ be a point of $Z\sb{\rep g}$.
We put 
$Q:=\pghg$,
which is  a point of $\tlZZ$.
Let $I\sb p$ denote the homogeneous ideal of $R$ defining the point $p$.
Then 
the dimension of the kernel of
$(d\alpha)\sb{Q}: T\sb Q\tlZZ \to T\sb{\alpha (Q)} \XX$ is equal to
\begin{equation}\label{eq:dimKer}
\dim \Fba-\dim\Ha-s+\dim (\Maa / (J\sb g \Maa+I\sb p h (\Mbb)))\sb 0.
\end{equation}
\end{sproposition}
\begin{sproof}
Since $\tlZZ$ is a closed subscheme of 
$\Hb \times \XX$
with $\rho\sprime\circ \pi$ and
$\alpha$ being the projections,
the kernel of $(d\alpha)\sb Q$ is mapped
isomorphically to a subspace of
$T\sb{\rep g}\Hb$
by the linear map $d (\rho\sprime \circ \pi)\sb Q$.
We will show that
this subspace
\begin{equation}\label{eq:subs1}
( d (\rho\sprime \circ \pi)\sb Q)\, (\Ker (d\alpha)\sb Q)
\;\;\subset \;\;
T\sb{\rep g} \Hb
\end{equation}
coincides
with
the subspace
\begin{equation}\label{eq:subs2}
(I\sb p \Mbb / J\sb{g} \Mbb)\sb 0  \;\cap\; 
\Ker  \rep h\sb g
\;\;\subset \;\; (\Mbb/J\sb g \Mbb)\sb 0
\end{equation}
under the identification~\eqref{eq:ZNZ}.
Let $g\sprime$ be an element of $\Mb$,
and  $\bar g\sprime$ 
the element $g\sprime  \mod (J\sb g \Mbb)\sb 0$ of $(\Mbb/J\sb g \Mbb)\sb 0$,
giving 
the corresponding displacement
${Z}\sb{\rep {g+\ve g\sprime}} \to\Spec \C [\ve]$
of $Z\sb{\rep g}$.
The subspace~\eqref{eq:subs1}
consists
of vectors
corresponding to 
infinitesimal displacements 
with  $p$ in ${Z}\sb{\rep {g+\ve g\sprime}}$
and with ${Z}\sb{\rep {g+\ve g\sprime}}$
remaining in $X\sb{\rep{h(g)}}$.
The displacement
${Z}\sb{\rep {g+\ve g\sprime}}$
contains $p$ 
if and only if $J\sb{g\sprime}\subset I\sb p$ holds,
which is equivalent to
$\bar g\sprime \in (I\sb p\Mbb/ J\sb g \Mbb)\sb 0$.
On the other hand,
by Proposition~\ref{prop:diagT},
the displacement ${Z}\sb{\rep {g+\ve g\sprime}}$
remains in $X\sb{\rep{h(g)}}$
if and only if
the corresponding vector of $T\sb{\rep g} \Hb$
is contained in  $\Ker \zeta$.
Since
$\zeta$ is identified with $ \rep h \sb g$,
this holds
if and only if
$\bar g\sprime\in \Ker \rep h \sb g$.
Therefore~\eqref{eq:subs1}
coincides
with~\eqref{eq:subs2}
by~\eqref{eq:ZNZ}.
The cokernel of the homomorphism
$$
(I\sb p \Mbb / J\sb{g} \Mbb)\sb 0  \;\hookrightarrow\;
(\Mbb/J\sb g \Mbb) \sb 0 \;\maprightsp{\rep h\sb g}\;
(\Maa/J\sb g \Maa)\sb 0
$$
is 
$(\Maa / (J\sb g\Maa +I\sb p h (\Mbb)))\sb 0$.
On the other hand, 
$\dim (\Mbb / I\sb p \Mbb)\sb 0$ is equal to $s$.
These lead us to the conclusion that
$\dim \Ker (d\alpha)\sb Q$ is equal to
$$
\dim (\Mbb / J\sb g \Mbb)\sb 0 -s - 
\dim(\Maa / J\sb g \Maa )\sb 0 +
\dim (\Maa / (J\sb g\Maa +I\sb p h (\Mbb)))\sb 0,
$$
which coincides with~\eqref{eq:dimKer} by Corollary~\ref{cor:dimFbaanddimCokerdrho} (1).
\end{sproof}
\par
%
%
%
%
In the sequel,
we use the following notation.
For positive integers $d$ and $e$, 
let $\Mat(d, e)$ denote the vector space
of all $d\times e$ matrices
with entries in $\C$,
and
$D(d, e)$ the Zariski closed subset of $\Mat (d, e)$
consisting of matrices
whose rank is less than $\min (d, e)$.
It is easy to see that
$D(d, e)$ is  irreducible.
We set
$$
o:=[1: 0: \dots: 0]\in \P\sp n.
$$
For a homogeneous polynomial $a\in R$,
we put
$$
a(o):=\hbox{the coefficient of $x\sb 0\sp{\deg a}$ in $a$}.
$$
Let $I\sb o$ be the homogeneous ideal of $R$
defining $o$ in $\P\sp n$:
$$
I\sb o:=\langle  x\sb1, \dots, x\sb n\rangle \subset R.
$$
We define linear maps
$\lambda\sb i : \IMa\to\C\sp n \;(i=1, \dots, r)$
and
$\mu\sb j : \IMb\to\C\sp n\; (j=1, \dots, s)$
by
$$
\lambda\sb i (f):=
\left(
\frac{\partial f\sb i}{\partial x\sb 1} (o),\dots,
\frac{\partial f\sb i}{\partial x\sb n} (o)
\right)
\quand
\mu\sb j (g):=
\left(
\frac{\partial g\sb j}{\partial x\sb 1} (o),\dots,
\frac{\partial g\sb j}{\partial x\sb n} (o)
\right).
$$
Let
$\lambda : \IMa \to \Mat (r, n)$
and 
$\mu : \IMb \to \Mat (s, n)$
be linear maps defined by 
$$
\lambda (f):=
\begin{pmatrix}
\lambda\sb 1 (f) \\
\vdots\\
\lambda \sb r (f)\\
\end{pmatrix}
\quand
\mu (g):=
\begin{pmatrix}
\mu\sb 1 (g) \\
\vdots\\ 
\mu \sb s (g)\end{pmatrix}.
$$
Both of $\lambda$ and $\mu$ are surjective.
We define
a linear map 
$\eta : \Hom (\Mbb, \Maa)\sb 0\to \Mat(r, s)$
by
$$
\eta (h) :=(h\sb{ij} (o)),
$$
when $h$ is expressed by an $r\times s$
matrix $(h\sb{ij})$ with $h\sb{ij}\in R\sb{a\sb i -b\sb j}$.
Note that,
if $g$ is an element of $(\Io \Mbb)\sb 0$,
then,
for any $h\in \HM$,
we have
$h(g)\in (I\sb o \Maa)\sb 0$
and
$\lambda (h (g)) =\eta (h) \cdot \mu (g)$.
\par
We define an $R$-submodule $N\sb{\aa}$ of $\Maa$
by
$$
N\sb{\aa} :=\bigoplus\sb{i=1}\sp{r-1} R(a\sb i) \oplus I\sb o (a\sb r).
$$
Note that
$\Ker \lambda\sb r=\IN$
holds in $\IMa$.
Note also that an element $h$ of $\HM$ is contained in $\HMN$ if and only if
the $r$-th row vector
of $\eta (h)$ is the zero vector.
We put
\begin{eqnarray*}
\IMbci &:=& (I\sb o \Mbb)\sb 0 \cap \Mbci, \\
\IMaci &:=& (I\sb o \Maa)\sb 0 \cap \Maci,  \\
\INci &:=& (I\sb o \Naa)\sb 0 \cap \Maci. 
\end{eqnarray*}
For $f\in (\Maa)\sb 0\spci$
and
$g\in (\Mbb)\sb 0\spci$,
we have the following:
\begin{equation}\label{eq:IoMaabb}
(o, \rep f)\in \XX \Longleftrightarrow
f\in (\Io\Maa)\sb 0\spci,
\;\;
(o, \rep g)\in \ZZ \Longleftrightarrow
g\in (\Io\Mbb)\sb 0\spci.\;\;
\end{equation}
Let $\Gamma$ be the Zariski closed subset
of $\XX$ consisting of critical points of $\phi :\XX\to\Ha$,
and  $\Gamma\sprime$
the Zariski closed subset
of $\ZZ$ consisting of critical points of $\phi\sprime :\ZZ\to\Hb$.
We put
$$
\Gamma\sb o := \tau\inv (o) \cap \Gamma,
\quad
\Gamma\sprime\sb o := \tau\sp{\prime\;-1} (o) \cap \Gamma\sprime.
$$
For $f\in (\Io\Maa)\sb 0\spci$
and
$g\in (\Io\Mbb)\sb 0\spci$,
we have the following:
\begin{equation}\label{eq:Gammao}
\begin{split}
&(o, \rep f)\in \Gamma\sb o \Longleftrightarrow
\lambda (f) \in D(r, n),\\
&(o, \rep g)\in \Gamma\sprime\sb o \Longleftrightarrow
\mu (g) \in D(s, n).
\end{split}
\end{equation}
If $f\in \INci$,  then $\lambda (f) \in D(r, n)$.
Hence
we can define a morphism
$\map{\gamma}{ \INci }{ \Gamma\sb o}$
by 
$$
\gamma (f):=(o, \rep f).
$$
\begin{sproposition}\label{prop:Gammao}
The Zariski closed subset $\Gamma\sb o$ of $\XX$ is
irreducible,
and the morphism $\gamma: \INci \to  \Gamma\sb o$ is dominant. 
\end{sproposition}
\begin{sproof}
By~\eqref{eq:IoMaabb} and~\eqref{eq:Gammao}, the map
$f\mapsto (o,\rep f)$
gives a surjective
morphism 
from $\lambda\inv (D(r, n))\cap \IMaci$
to $\Gamma\sb o$.
Because $\lambda$ is a surjective 
linear map
and $D(r, n)$ is irreducible,
$\lambda\inv (D(r, n))$ is also
irreducible.
Since $\lambda\inv (D(r, n))\cap \IMaci$ is
Zariski open in  $\lambda\inv (D(r, n))$, 
$\Gamma\sb o$ is also irreducible.
Let $f$ be a general element
of $\lambda\inv (D(r, n))$.
Then
$\lambda (f)$ is of rank $r-1$,
and the vector $\lambda\sb r (f)$
can be written 
as a linear combination
of $\lambda\sb{1} (f), \dots, \lambda\sb{r-1} (f)$.
Since
$a\sb r\ge a\sb i$ 
for $i<r$,
there exist
 homogeneous polynomials
$c\sb 1, \dots, c\sb{r-1}$
with $c\sb i\in R\sb{a\sb r-a\sb i}$
such that,
if we put
$$
f\sb r\sprime:=f\sb r -c\sb1 f\sb1 -\cdots -c\sb{r-1} f\sb{r-1}
\quand
f\sprime:=\transpose{(f\sb 1, \dots, f\sb{r-1}, f\sb r\sprime)},
$$
then
$\lambda\sb r (f\sprime)=0$ holds,
which means $f\sprime \in \IN$.
From  $J\sb f=J\sb{f\sprime}$,
we conclude that $(o, \rep f)=(o, \rep{f\sprime})$ belongs to the image of $\gamma$.
Since $(o, \rep{f})$ is a general point of $\Gamma\sb o$,
the morphism $\gamma$ is dominant.
\end{sproof}
\begin{sremark}
The irreducibility of $\Gamma$ and
that of $\Sa=\phi (\Gamma)$ follow
from  Proposition~\ref{prop:Gammao}
and  the action of $\PGL (n+1)$ on the diagram~\eqref{diag:univ}.
\end{sremark}
\begin{scorollary}\label{cor:generalinGammao}
Suppose that  $(o,\rep f)$ is a general point of $\Gamma\sb o$.
Then the singular locus of $X\sb{\rep f}$ consists of only one  point $o$,
which is a hypersurface singularity 
of $X\sb{\rep f}$ with non-degenerate Hessian.
\noproof
\end{scorollary}
We put
$$
\Xi 
:=\alpha\inv (\Gamma)\setminus 
(\alpha\inv (\Gamma)\cap \beta\inv (\Gamma\sprime))
\quand
\Xi\sb o :=(\tau\circ\alpha)\inv (o) \cap \Xi,
$$
which are locally closed subsets of $\tlZZ$ (possibly empty).
A point $(o, \rep g, \rep f)$ of $(\tau\circ\alpha)\inv (o)\subset \tlZZ$
is contained in $\Xi\sb o$ if and only if 
$X\sb{\rep f}$ is singular at $o$  and 
$Z\sb{\rep g}$ is smooth at $o$.
The morphism
$\alpha : \tlZZ\to\XX$ induces a morphism 
$\alpha|\Xi\sb o : \Xi\sb o \to \Gamma\sb o$.
\begin{sremark}
Invoking  the action of $\PGL (n+1)$
on the diagram~\eqref{diag:univ},
we can paraphrase
the second condition of Main Theorem
into the condition  that
$\alpha|\Xi\sb o : \Xi\sb o \to \Gamma\sb o$
is dominant.
\end{sremark}
We define  a linear subspace $V$ of
$U=(\Mbb)\sb 0 \times \Hom (\Mbb, \Maa)\sb 0$
by
$$
V :=(I\sb o \Mbb)\sb 0 \times \Hom (\Mbb, N\sb{\aa})\sb 0.
$$
We then put $V\spci := V\cap U\spci$ and 
$$
V\spdi :=\{\,(g, h)\in V\spci  \mid \mu(g)\notin D(s, n) \}=
 \{\,(g, h)\in V\spci  \mid \textrm{$Z\sb{\rep{g}}$ is smooth
at $o$} \}.
$$
By definition,
$V\spdi$ is Zariski open in the vector space $V$,
but may possibly be empty.
Recall that
$\nu : U \to \Ma$
is the morphism defined by $\nu (g, h)=h(g)$.
We have a morphism
$$
\nu|V:  V \to \IN
\quand
\nu|V\spdi:  V\spdi \to \INci,
$$
which  are the restrictions of  $\nu$ to $V$ and $V\spdi$, respectively.
By definition again,
if $(g, h)\in V\spdi$,
then
$\oghg\in \Xi\sb o$.
Let  $\xi : V\spdi \to \Xi\sb o$ be the morphism defined by
$$
\xi (g, h) :=\oghg.
$$
Then we obtain the following commutative diagram:
\begin{equation}\label{diag:Vd}
\begin{array}{ccc}
\;\;V\spdi\;\; &\maprightsp{\nu|V\spdi} &\; \INci \\
 \mapdownleft{\xi} & & \mapdownright{\gamma}\\
\Xi\sb o & \maprightsb{\alpha|\Xi\sb o} & \Gamma\sb o.\\
\end{array}
\end{equation}
\begin{sproposition}\label{prop:alphaXio}
The morphism
$\alpha |\Xi\sb o : \Xi\sb o \to \Gamma\sb o$
is dominant if and only if
$\nu |V\spdi : V\spdi\to \INci$ is dominant.
\end{sproposition}
\begin{sproof}
Since $\gamma$ is dominant by Proposition~\ref{prop:Gammao},
the commutativity of the digram~\eqref{diag:Vd}
implies that,
if $\nu|V\spdi$ is dominant,
then so is $\alpha | \Xi\sb o$.
Suppose conversely
that $\alpha|\Xi\sb o$
is dominant.
Let $f$ be a general point of
$\INci$.
Since $\gamma$ is dominant,
$(o, \rep f)$
is a general point of $\Gamma\sbo$,
and hence $(o, \rep f)$ is in the image of $\alpha |\Xi\sb o$.
Thus
there exists an element
$g\in \IMbci$
such that
$(o, \rep g, \rep f) \in \Xi\sb o$,
which implies that
$\mu (g)$ is {\em not} contained in $ D(s, n)$, 
and that  there exists an element 
$h\in \HM$ that satisfies $h(g)=f$.
From $\lambda\sb r(f)=0$ and $\eta (h)\cdot \mu(g)=\lambda (f)$,
the linear independence
of the row vectors of $\mu (g)$
implies that the $r$-th row vector
of $\eta (h)$ is a zero vector.
Therefore $h$ is in fact
an element of $\HMN$,
which means $(g, h)\in V\spdi$.
Hence the general point 
$f=h(g)$ of $\INci$ is contained in the image of $\nu|V\spdi$.
\end{sproof}
\begin{sproposition}\label{prop:uniqueirred}
Suppose that $\alpha |\Xi\sb o : \Xi\sb o \to \Gamma\sb o$
is dominant.
Then there exists a unique
irreducible component
$\Xi\sb o\sprime$ of $\Xi\sb o$
such that the restriction 
$\alpha |\Xi\sb o\sprime : \Xi\sb o\sprime \to \Gamma\sb o$
of $\alpha |\Xi\sb o$ to $\Xi\sb o\sprime$
is dominant.
The closure of the image of $\xi: V\spdi \to \Xi\sb o$ in $\Xi\sb o$
coincides with $\Xi\sb o\sprime$.
\end{sproposition}
\begin{sproof}
Since $\Gamma\sb o$ is irreducible,
there exists at least one
irreducible component $\Xi\sb o\sprime$ of $\Xi\sb o$
that is mapped dominantly onto $\Gamma\sb o$ by $\alpha |\Xi\sb o$.
Let $(o, \rep g, \rep f)$ be a general point of $\Xi\sb o\sprime$.
Then $\alpha(o, \rep{g}, \rep{f})=(o, \rep f)$ is a general point of $\Gamma\sb o$.
Since $\gamma$ is dominant,
we can assume
that $(o, \rep{f})$ is in the image of $\gamma$;
that is, $f$ is an element of $\INci$.
Let $h\in \HM$ be a homomorphism satisfying  $h(g)=f$.
From $\mu (g)\notin D(s, n)$ and $\lambda\sb r (f)=0$,
we see that 
$h$ actually is 
an element of $\HMN$.
Hence $(g, h)$ is a point of $V\spdi$,
which is  mapped to the general point $(o, \rep g, \rep f)$ of $\Xi\sb o\sprime$ by $\xi$.
Therefore $\Xi\sb o\sprime$ is the closure
of the image of $\xi: V\spdi \to \Xi\sb o$ in $\Xi\sb o$.
Since $V\spdi$ is irreducible,
the uniqueness of $\Xi\sb o\sprime$, as well as the second assertion,  is proved.
\end{sproof}
For an element $(g, h)$ of $U$,
we define a linear map
$\delta\sbgh : U \to \Ma$ by
$$
\delta\sbgh (G, H):= H(g)+ h(G)
\quad (\, G\in \Mb,\; H\in\HM\,).
$$
Under the natural isomorphisms
$T\sbgh U\cong U$ and $T\sb{\nu(g, h)}\Ma\cong \Ma$,
the linear map $\delta\sbgh$ is equal to 
$$
(d\nu)\sbgh : T\sbgh U \to T\sb{\nu(g, h)}\Ma.
$$
By definition, we  have
\begin{eqnarray}\label{eq:imdelta1}
\delta\sbgh (U) \;&=&\; (J\sb g \Maa + h (\Mbb))\sb 0, \\
\label{eq:imdelta3}
\delta\sbgh (V) \;&=&\; (J\sb g N\sb{\aa} +I\sb o h (\Mbb))\sb 0,
\quand \\
\label{eq:imdeltaV}
(g, h)\in V\;&\Longrightarrow&\; \,\delta\sbgh (U) \subseteq \Na,\;\;\delta\sbgh (V) \subseteq \IN.
\end{eqnarray}
\begin{sproposition}\label{prop:equivconds}
Suppose  $a\sb r \ge b\sb s$.
Then 
the following conditions on $(n, \aa, \bb)$ are
equivalent to each other\itcolon
\par
{\rm (i)}
				The morphism $\alpha |\Xi\sb o : \Xi\sb o \to \Gamma\sb o$ is dominant.
\par
{\rm (ii)}
				If $(g, h)\in V$ is general, then $\delta\sbgh (V)$ coincides with $\IN$.
\par
{\rm (iii)}
    If $(g, h)\in V$ is general, then the following holds\itcolon
    $$
       \dim (\Maa / (J\sb g \Maa + I\sb o h (\Mbb)))\sb 0 =n+r-s.
    $$
\par
{\rm (iv)}
				There exists at least one $(g, h)\in V$ such that  
    $$
      \dim (\Maa / (J\sb g \Maa + I\sb o h (\Mbb)))\sb 0 \le n+r-s.
    $$
\end{sproposition}
\begin{sproof}
First we show the following:
\par
\smallskip
\noindent
{\bf Claim}
(1)
For any $(g, h) \in V$,
$\dim (\Maa/ (J\sb g \Maa+\Io h (\Mbb)))\sb 0$ 
is larger than or equal to $n+r-s$.
(2)
If $(g, h)\in V$ is chosen generally, then
$\dim (J\sb g \Maa+\Io h (\Mbb))\sb 0 $ is equal to
$\dim (J\sb g \Naa+\Io h (\Mbb))\sb 0+s$.
\par
\smallskip
\noindent
Let $(g, h)$ be an arbitrary element of $V$.
Then 
$(\Io h (\Mbb))\sb 0 $ is contained in $ (\Io\Naa)\sb 0 =\Ker\lambda\sb r$.
On the other hand,
if $f\in (J\sb g \Maa)\sb 0$,
then the $r$-th component $f\sb r$ of $f$ is written
as $g\sb 1 k\sb 1 + \cdots + g\sb s k\sb s$
%
%
%
%
with $k\sb j \in R\sb{a\sb r -b\sb j}$,
and $\lambda\sb r (f)$
is equal to
\begin{equation}\label{eq:kg}
k\sb 1 (o) \mu\sb 1 (g) + \cdots + k\sb s (o) \mu\sb s (g).
\end{equation}
Hence the image of
$(J\sb g \Maa+\Io h (\Mbb))\sb 0$
by $\lambda\sb r$
is spanned by
$\mu\sb 1 (g)$, \dots, $\mu\sb s (g)$,
and therefore is of dimension $\le s$.
On the other hand,
$\Ker\lambda\sb r = (\Io\Naa)\sb 0$ is of codimension $n+r$ in $\Ma$.
Hence we obtain
$$
\dim (J\sb g \Maa+\Io h (\Mbb))\sb 0
\le
\dim \Ker \lambda\sb r + s
=
\dim \Ma -n-r+s,
$$
which implies Claim~(1).
\par
Let $(g, h)$ be a general element of $V$.
Because $g$ is general in $\IMb$,
the vectors 
$\mu\sb 1 (g), \dots, \mu\sb s (g)$
are linearly independent.
Let $f$ be an element of $(J\sb g \Maa)\sb 0$.
By the assumption $a\sb r\ge b\sb s$,
the degrees $a\sb r-b\sb j$ of the polynomials $k\sb j$ in the expression
$f\sb r=g\sb 1 k\sb 1 + \cdots + g\sb s k\sb s$
are non-negative for all $j\le s$.
Therefore the coefficients
$k\sb j (o)$ in~\eqref{eq:kg} can take any values.
Hence
the image
of
$(J\sb g \Maa+\Io h (\Mbb))\sb 0$
by $\lambda\sb r$
is of dimension exactly $s$.
Moreover,
if  $f\in \Ker\lambda\sb r$,
then
$k\sb 1 (o)=\cdots = k\sb s (o)=0$
holds.
Hence we have 
$\Ker\lambda\sb r \subseteq (J\sb g \Naa)\sb 0$.
Because
$(\Io h (\Mbb))\sb 0 \subseteq \Ker\lambda\sb r$, 
we have
$$
(J\sb g \Maa+\Io h (\Mbb))\sb 0\cap \Ker\lambda\sb r
=(J\sb g \Naa+\Io h (\Mbb))\sb 0.
$$
Therefore Claim~(2) is proved.
\par
Since
$\dim (\Maa/ (J\sb g \Maa+\Io h (\Mbb)))\sb 0$ 
is 
an upper semi-continuous function of $(g, h)$,
Claim~(1) implies that the conditions
(iii) and (iv) are equivalent.
By~\eqref{eq:imdelta3} and~\eqref{eq:imdeltaV}, 
the following inequality holds for any $(g, h) \in V$:
\begin{equation}\label{eq:ineq}
\begin{split}
\dim \Ma / \delta\sb{(g,h)} (V)&=\dim (\Maa / (J\sb g \Naa+\Io h (\Mbb)))\sb 0 \\
 & \ge
\dim (\Maa/ \Io\Naa)\sb 0 =n+r.
\end{split}
\end{equation}
The condition (ii) is satisfied if and only if
the {\em equality} in~\eqref{eq:ineq} holds
for a general $(g, h)\in V$.
The equivalence of 
the conditions (ii) and (iii)
now follows
from Claim~(2).
\par
By Proposition~\ref{prop:alphaXio},
the condition (i)
is equivalent to the following:
\par
(i)${}\sprime$
The morphism $\nu | V\spdi : V\spdi \to \INci$ is dominant.
\par\noindent
On the other hand,
since
$\delta\sbgh$
is equal to
$(d\nu)\sbgh : T\sbgh U \to T\sb{\nu(g, h)}\Ma$
via the natural identifications
$T\sbgh U\cong U$ and $T\sb{\nu(g, h)}\Ma\cong \Ma$,
the condition (ii)
is equivalent to the following:
\par
(ii)${}\sprime$
The morphism $\nu | V : V \to \IN$ is dominant.
\par\noindent
Since $\INci$ is Zariski open dense in $\IN$,
the implication
$\hbox{\rm (i)}\Rightarrow \hbox{\rm (ii)}$
is obvious.
Since $V\spdi$ is Zariski open in $V$,
the implication 
$\hbox{\rm (ii)}\Rightarrow \hbox{\rm (i)}$
follows if we show  that
$V\spdi$ is non-empty
under  the condition (ii).
Suppose that
the condition (ii) is fulfilled.
Let $(g, h)$ be a general element of $V$.
Since $g$ is general in $\IMb$,
the ideal $J\sb g$ defines a complete intersection of multi-degree $\bb$
passing through $o$,
and $\mu (g)$ is of rank $s$.
By (ii)${}\sprime$,
$h(g)$  is a general element of $\IN$,
and hence
$J\sb{h(g)}$ defines a complete intersection of multi-degree $\aa$ passing through $o$ and singular at $o$.
Thus we have $(g, h)\in V\spdi$.
\end{sproof}
\section{Proof of Main Theorem}\label{sec:proofMT}
First we prepare two easy lemmas.
\par
Let $L\sb 1$ and $L\sb 2$ be finite-dimensional  vector spaces,
and  let $\Hom (L\sb 1, L\sb 2)$  be the vector space of
linear maps from $L\sb 1$ to $L\sb 2$.
For $\varphi \in \Hom(L\sb 1, L\sb 2)$,
we have a canonical identification
\begin{equation}\label{eq:can}
T\sb{\varphi} \Hom (L\sb 1, L\sb 2) \cong \Hom (L\sb 1, L\sb 2).
\end{equation}
Let $S\sb k$ be the closed subscheme
of $\Hom(L\sb 1, L\sb 2)$
defined as  common zeros 
of all $(k+1)$-minors
of the matrices expressing the linear maps
in terms of certain bases of $L\sb 1$ and $L\sb 2$.
\begin{slemma}\label{lem:lin}
Let $\varphi\sb 0$ be a point of $S\sb k\setminus S\sb{k-1}$.
An element $\varphi$ of $\Hom (L\sb 1, L\sb 2)$
is contained in the subspace 
$T\sb{\varphi\sb 0} S\sb k$ of $ T\sb{\varphi\sb 0} \Hom (L\sb 1, L\sb 2)$
under the identification~\eqref{eq:can}
if and only if $\varphi (\Ker \varphi\sb 0)$ is contained in $\Im \varphi\sb 0$.
\end{slemma}
\begin{sproof}
We can choose bases of $L\sb 1$ and $L\sb 2$ in such a way that
$\varphi\sb 0$
is expressed by the matrix
{\small$
\left(
\begin{array}{@{\,}c|c@{\,}}
I\sb k & O \\ \hline
O & O
\end{array}
\right).
$
}
Suppose that $\varphi$ is expressed by the matrix
{\small
$
\left(
\begin{array}{@{\,}c|c@{\,}}
A & B \\ \hline
C & D
\end{array}
\right)
$
}
under these bases.
Then 
$\varphi$
is contained in
$T\sb{\varphi\sb 0} S\sb k$
under the identification~\eqref{eq:can}
if and only if the matrix
{\small
$
\left(
\begin{array}{@{\,}c|c@{\,}}
I\sb k +\ve A & \ve B \\ \hline
\ve C & \ve D
\end{array}
\right) 
$
}
is of rank $k$,
where $\ve$ is the dual number; $\ve\sp 2=0$.
This matrix is of rank $k$ if and only if $D=0$,
which is equivalent to $\varphi (\Ker \varphi\sb 0) \subseteq \Im \varphi\sb 0$.
\end{sproof}
\par
Let $X$ and $Y$ be connected complex manifolds,
$Z$ an irreducible locally closed  analytic
subspace of $Y$,
$\psi : X\to Y$ a holomorphic map,
and 
$p$ a point of $\psi\inv (Z)$.
\begin{slemma}\label{lem:cut}
Suppose that $Z$ is smooth at $\psi (p)$,
and that we have
\begin{equation}\label{eq:trans}
T\sb{\psi (p)} Z \cap \Im (d\psi)\sb p =0
\quand
T\sb{\psi (p)} Z + \Im (d\psi)\sb p =T\sb{\psi(p)} Y.
\end{equation}
Then $\psi\inv (Z)$ is smooth at $p$.
Moreover,
the dimension of $\psi\inv (Z)$ at $p$
is equal to $\dim X - \dim Y +\dim Z$.
\end{slemma}
\begin{sproof}
By~\eqref{eq:trans},
we have
$T\sb p \psi\inv (Z) = (d\psi)\sb p\inv (T\sb{\psi (p)} Z) =\Ker (d\psi)\sb p$,
and hence 
$\dim T\sb p \psi\inv (Z) $ is equal to $ \dim T\sb p X -\dim \Im  (d\psi)\sb p$,
which is then 
equal to 
$\dim X-  \dim Y +\dim Z$.
On the other hand,
the codimension of $\psi\inv (Z)$ in $X$ at $p$
is less than or equal to
the codimension of $Z$ in $Y$ at $\psi (p)$.
Combining these facts,
we get the hoped-for results.
\end{sproof}
\par
From now on, we assume that $(n, \aa, \bb)$ satisfies the conditions required in  Main Theorem.
In particular, the morphism $\alpha |\Xi\sb o : \Xi\sb o \to \Gamma\sb o$
is dominant.
Let $\Xi\sb o\sprime$
be the unique irreducible component  of $\Xi\sb o$
that is mapped dominantly onto $\Gamma\sb o$ by $\alpha |\Xi\sb o$
(Proposition~\ref{prop:uniqueirred}).
\begin{sproposition}\label{prop:dims2}
Let $Q=(o, \rep g, \rep f)$
be a general point of $\Xi\sb o\sprime$.
Then the following hold\itcolon
\par
{\rm (1)}
The morphism $\rho : \Fba \to \Ha$ is dominant.
\par
{\rm (2)}
The kernel of
$(d\alpha)\sb Q : T\sb Q\tlZZ \to T\sb{\alpha (Q)} \XX$
is of dimension equal to 
$$
\dim \Fba -\dim\Ha -m+2l.
$$
%
%
%
%
\par
{\rm (3)}
The image of 
$(d\rho)\sb{\pi (Q)} : T\sb{\pi (Q)} \Fba \to T\sb{\rep f} \Ha$
is of codimension  $1$.
\end{sproposition}
\begin{sproof}
First of all, note that
$\Fba$ is non-empty because $\Xi\sb o$ is non-empty.
Note also that $ V\spdi$ is non-empty by Proposition~\ref{prop:alphaXio},
and hence is Zariski open dense in $V$.
By Proposition~\ref{prop:uniqueirred},
the general point $Q$ of $\Xi\sb o\sprime$ is the image of a general point of 
$ V\spdi$
by $\xi : V\spdi \to \Xi\sb o$.
Therefore
we can choose a general point $(g, h)$ of $ V$
first, 
and then put $Q:=\xi (g, h)=\oghg$.
\par
%
%
We  start with the proof of (3).
By Proposition~\ref{prop:equivconds},
we have
\begin{equation}\label{eq:ImIN}
\delta\sb{(g, h)} (V) =(J\sb g \Naa+\Io h (\Mbb))\sb 0 =\IN.
\end{equation}
In particular,
we have
\begin{equation}\label{eq:JgIo}
\dim (\Maa/ (J\sb g \Naa+\Io h (\Mbb)))\sb 0 =n+r.
\end{equation}
By Corollary~\ref{cor:dimFbaanddimCokerdrho} (2),
to arrive at
$\dim \Coker (d\rho)\sb{\pi (Q)}=1$,
all we have to show is 
$$
\dim (\Maa/ (J\sb g \Maa+\ h (\Mbb)))\sb 0 =1,
$$
which is equivalent to
\begin{equation}\label{eq:cokerdim}
\dim \delta\sb{(g, h)} (U) / \delta\sb{(g, h)} (V)
=n+r-1,
\end{equation}
because of~\eqref{eq:imdelta1}, \eqref{eq:imdelta3} and~\eqref{eq:JgIo}.
We define  a linear map
$\tilde\lambda\sb r : \Ma \to \C\sp r\times \C\sp n$
by
$$
\tilde\lambda\sb r (f)
:=
\Bigl ( (f\sb 1 (o), \dots, f\sb r (o)), 
\Bigl(\frac{\partial f\sb r}{\partial x\sb 1} (o), \dots, \frac{\partial f\sb r}{\partial x\sb n} (o)\Bigr)
\Bigr).
$$
Then 
$\delta\sb{(g, h)} (V) =\IN=\Ker \tilde\lambda\sb r$ holds from~\eqref{eq:ImIN}.
Moreover,
we have
$\delta\sb{(g, h)} (U)\subseteq\Na$ by~\eqref{eq:imdeltaV}, and $\dim\tilde\lambda\sb r (\Na)=n+r-1$.
Therefore~\eqref{eq:cokerdim}
and the following two conditions are 
equivalent to each other:
\begin{eqnarray}
\label{eq:N}
&\delta\sb{(g, h)} (U) = \Na,
\\
\label{eq:cokerdim2}
&\dim \tilde\lambda\sb r (\delta\sb{(g, h)} (U))  \ge n+r-1.
\end{eqnarray}
We will prove~\eqref{eq:cokerdim} by showing that the inequality~\eqref{eq:cokerdim2}
holds.
For $\nu=1, \dots, s$,
we define $\gamma\sp{(\nu)}=\transpose{(\gamma\sp{(\nu)}\sb1, \dots, \gamma\sp{(\nu)}\sb s)} \in \Mb$ 
and
$\eta\sp{(\nu)}=(\eta\sp{(\nu)}\sb{ij})\in \HM$ by
$$
\gamma\sp{(\nu)}\sb j:=
\begin{cases}
0 &\textrm{if $j \ne \nu$}\\
x\sb 0\sp{b\sb j} & \textrm{if $j=\nu  $} \\
\end{cases}
\quand
\eta\sb{ij}\sp{(\nu)}
:=
\begin{cases}
0 & \textrm{ if $(i, j) \ne (r, \nu)$} \\
x\sb 0 \sp{a\sb r -b\sb{\nu}} & \textrm{ if $(i, j) = (r, \nu)$}.
\end{cases}
$$
Note that $a\sb r \ge b\sb\nu$ by~\eqref{eq:num1}.
We then define
$v\sp{(\nu)}, w\sp{(\nu)}\in \delta\sb{(g, h)} (U)$
by
\begin{equation*}
\begin{split}
v\sp{(\nu)} &:=\delta\sb{(g, h)} (\gamma\sp{(\nu)}, 0)= h (\gamma\sp{(\nu)}),\\
w\sp{(\nu)} &:= \delta\sb{(g, h)} (0, \eta\sp{(\nu)})
=\eta\sp{(\nu)} (g) =
\transpose{(0, \dots, 0, g\sb{\nu} x\sb 0\sp{a\sb r-b\sb{\nu}})}. \\
\end{split}
\end{equation*}
Then we have
\begin{equation}\label{eq:lambdavw}
\begin{split}
\tilde\lambda\sb r (v\sp{(\nu)})  &=
\Bigl(
(h\sb{1\nu} (o), \dots, h\sb{r\nu} (o)), 
\Bigl(
\frac{\partial h\sb{r \nu}}{\partial x\sb 1} (o),
\dots, 
\frac{\partial h\sb{r \nu}}{\partial x\sb n} (o)
\Bigr)
\Bigr),
\\
\tilde\lambda\sb r (w\sp{(\nu)})  &=
\Bigl(
(0, \dots, 0), 
\Bigl(
\frac{\partial g\sb\nu}{\partial x\sb 1 } (o),
\dots,
\frac{\partial g\sb\nu}{\partial x\sb n } (o)
\Bigr)\Bigr).
\end{split}
\end{equation}
In order to prove~\eqref{eq:cokerdim2},
it is enough to show 
that the vectors
$\tilde\lambda\sb r (v\sp{(\nu)}) $ and 
$\tilde\lambda\sb r (w\sp{(\nu)})$
span a hyperplane in $\C\sp r\times\C\sp n$.
Since $(g, h)$ is general in $V$, 
the coefficients
$h\sb{i\nu} (o) $, ${\partial h\sb{r\nu}}/{\partial x\sb j} (o)$ and ${\partial g\sb\nu}/{\partial x\sb j } (o)$
of the homogeneous polynomials $g\sb \nu$ and $h\sb{i\nu}$
that appear in~\eqref{eq:lambdavw}
are general except for the following
restrictions:
\begin{equation*}\label{eq:numres}
\begin{split}
&h\sb{r\nu} (o)=0\;\;(1\le \nu \le s) , \qquad h\sb{i\nu} (o)=0 \;\;\hbox{\textrm{if $a\sb i < b\sb \nu$}} \quand\\
&\frac{\partial h\sb{r\nu}}{\partial x\sb j} (o)=0 \;\;(1\le j \le n)\quad\hbox{\textrm{if $a\sb r=b\sb {\nu}$}}.
\end{split}
\end{equation*}
%
%
%
%
%
%
%
%
%
Let $\Lambda$ be 
the $2s \times (n+r)$ matrix  
whose row vectors are  $\tilde\lambda\sb r (v\sp{(\nu)}) $ and 
$\tilde\lambda\sb r (w\sp{(\nu)}) $.
Then $\Lambda$
is of the shape depicted in Figure~\ref{fig:matrix},
in which the entries in the submatrices marked with $*$
are general,
and the $(\nu, i)$-component
in the submatrix marked with $\sharp$
is general except for the restriction
that it must be zero if $a\sb i<b\sb {\nu}$.
\begin{figure}[t]
\includegraphics[height=10.3cm, keepaspectratio]{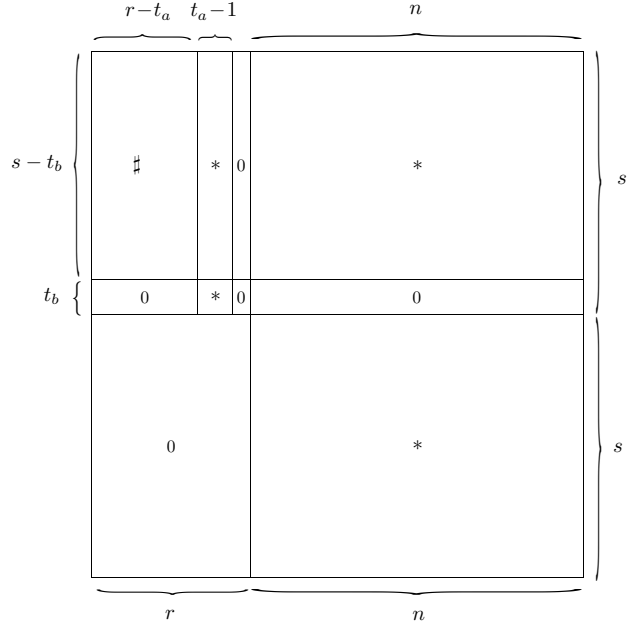}
\vskip -1.5cm
\caption{The shape of a $2s\times (n+r)$ matrix}
\label{fig:matrix}
\end{figure}
Since the rank of a matrix is a lower semi-continuous
function of entries,
in order to prove that $\Lambda$ is of rank $n+r-1$,
it is enough to show that
there exists at least one matrix of rank $n+r-1$
with the  shape Figure~\ref{fig:matrix}.
The  condition~\eqref{eq:num1}
implies  that $r-t\sb a \le s-t\sb b$ holds,
and that the $(i, i)$-component   of the submatrix $\sharp$
is subject to no restrictions for $i=1, \dots, r-t\sb a$.
The  condition~\eqref{eq:num2}
implies $n+r< 2 s$ and $n+r+t\sb b -t\sb a \le 2 s$.
Therefore we can define
a $2s \times (n+r)$ matrix $C$ of the shape Figure~\ref{fig:matrix}
by Table~\ref{tab:defC},
where $c\sb i$ is the $i$-th column vector of $C$ and
$e\sb {\mu}$ is the column vector of dimension $2s$
whose $\nu$-th component is $\delta\sb{\mu\nu}$
(Kronecker's delta symbol).
\begin{table}[b]
\fbox{
\begin{minipage}{54mm}
{\small
When $t\sb a-1 \ge t\sb b$,
$$
c\sb i:=\begin{cases}
e\sb i & \textrm{if $1\le i \le r-t\sb a$} \\
e\sb{s-r+1+i} & \textrm{if $r-t\sb a < i <r$} \\
0 & \textrm{if $i=r$} \\
e\sb{i-t\sb a}& \textrm{if $r<i\le s+1$} \\
e\sb{i-1}& \textrm{if $s+1 < i \le n+r$.} \\
\end{cases}
$$
}
\end{minipage}
}
\fbox{
\begin{minipage}{61mm}
{\small
When $t\sb a-1 < t\sb b$,
$$
c\sb i:=\begin{cases}
e\sb i & \textrm{if $1\le i \le r-t\sb a$} \\
e\sb{s-r+1+i} & \textrm{if $r-t\sb a < i <r$} \\
0 & \textrm{if $i=r$} \\
e\sb{i-t\sb a}& \textrm{if $r<i\le s+t\sb a -t\sb b$} \\
e\sb{i+t\sb b-t\sb a}& \textrm{if $s+t\sb a - t\sb b < i \le n+r$.} \\
\end{cases}
$$
}
\end{minipage}
}
\phantom{aaa}
\caption{Definition of $C$}
\label{tab:defC}
\end{table}
It is easy to see that $C$ is of rank $n+r-1$.
Hence~\eqref{eq:cokerdim2}, and also~\eqref{eq:cokerdim} and~\eqref{eq:N}, are proved.
\par
Next we  prove (1).
Since both of $\Fba$ and $\Ha$ are  smooth  and irreducible,
it is enough  to show that
the morphism $\rho: \Fba \to\Ha$ is a submersion
at a general point of $\Fba$.
By Corollary~\ref{cor:dimFbaanddimCokerdrho} (2),
it is therefore enough to prove that 
the following equality holds
for a general $(\tilde g, \tilde h) \in U$:
\begin{equation}\label{eq:rhosurj}
\dim (\Maa / ( J\sb{\tilde g }\Maa + \tilde h (\Mbb)))\sb 0 =0.
\end{equation}
Since the left-hand side of~\eqref{eq:rhosurj}
is an upper semi-continuous function of
$(\tilde g,\tilde  h)\in U$, 
it suffices to show that 
there exists at least one $(\tilde g,\tilde  h)\in U$
for which~\eqref{eq:rhosurj} holds.
We will find $(\tilde g, \tilde h)$ satisfying~\eqref{eq:rhosurj}
in a small neighborhood of the chosen point $(g, h)$ in $U$.
From~\eqref{eq:N},
we have
$$
\dim (\Maa / ( J\sb{\tilde  g} \Maa + \tilde h (\Mbb)))\sb 0 \le \dim (\Maa/\Naa)\sb 0= 1
$$
for any $(\tilde g,\tilde h)$ in a small neighborhood of $(g, h)$ in $U$.
We suppose the following:
\begin{equation}\label{eq:contradiction}
\begin{split}
\textrm{
 $\delta\sb{(\tilde g, \tilde h)} (U) =( J\sb{\tilde g} \Maa + \tilde h (\Mbb))\sb 0$
is of codimension $1$ in $\Ma$} \\
\textrm{
for any $(\tilde g, \tilde h)$
in a small neighborhood of $(g, h)$ in $U$,
}
\end{split}
\end{equation}
and will derive a contradiction.
For a sequence $c=(c\sb 1, \dots, c\sb s)$
of complex numbers,
we define $\eta\sp c =(\eta\sp c\sb{ij} )\in \HM$
by
$$
\eta\sp c\sb{ij} :=\begin{cases}
0 & \textrm{if $i<r$} \\
c\sb j x\sb 0\sp{a\sb r-b\sb j}& \textrm{if $i=r$}, \\
\end{cases}
$$
and consider the infinitesimal deformation
$(g, h)+\ve (0, \eta\sp c)$ of
$(g, h)$ in $U$,
where $\ve$ is the dual number.
Here we use the  condition
$a\sb r\ge b\sb j$ again.
By~\eqref{eq:N}, Lemma~\ref{lem:lin} and the assumption~\eqref{eq:contradiction},
we have
$$
\delta\sb{(0, \eta\sp c)} (\Ker \delta\sb{(g, h)})\subseteq \Im \delta\sb{(g, h)} =\Na
$$
for any $c$,
which means
that,
if $(G, H)\in \Ker \delta\sb{(g, h)}$,
then 
$\eta\sp c (G)\in \Na$ for any $c$.
Hence we have
\begin{equation}\label{eq:implication}
(G, H)\in \Ker \delta\sb{(g, h)}
\quad\Longrightarrow\quad
G\in \IMb.
\end{equation}
Because~\eqref{eq:num2} implies 
$2s> n+r$,
there exists a non-trivial linear relation
$$
\sum\sb{\nu=1}\sp s \alpha\sb{\nu} \tilde\lambda \sb r (v\sp{(\nu)})
+
\sum\sb{\nu=1}\sp s \beta\sb{\nu} \tilde\lambda \sb r (w\sp{(\nu)})
=0
\qquad (\alpha\sb\nu, \beta\sb\nu\in \C)
$$
among the vectors~\eqref{eq:lambdavw} in $\C\sp r \times \C\sp n$. 
Since $g$ is general in $\IMb$
and $s<n$,
the vectors
$\tilde\lambda \sb r (w\sp{(\nu)})$ $(\nu=1, \dots, s)$
are linearly independent,
and hence at least one of 
$\alpha\sb 1, \dots, \alpha\sb s$ is non-zero.
We put
$$
(G\sb 1, H\sb 1):=
\left(\sum\sb{\nu=1}\sp s \alpha\sb{\nu} \gamma\sp{(\nu)}, 
\sum\sb{\nu=1}\sp s \beta\sb{\nu} \eta\sp{(\nu)}\right) \in U.
$$             
Then we have
$$
\delta\sb{(g, h)} (G\sb 1, H\sb 1)=
\sum\sb{\nu=1}\sp{s} \alpha\sb\nu v\sp{(\nu)}+
\sum\sb{\nu=1}\sp{s} \beta\sb\nu w\sp{(\nu)}
\in \Ker\tilde\lambda\sb r =\IN =\delta\sb{(g, h)}
(V), 
$$
where the last equality follows from~\eqref{eq:ImIN}.
Hence there exists
$(G\sb 2, H\sb 2)\in V$
such that $
(G\sb 1-G\sb 2,H\sb 1-H\sb 2)\in\Ker \delta\sb{(g, h)}$.
On the other hand,
since $G\sb 2 \in \IMb$
and 
at least one of 
$\alpha\sb 1, \dots, \alpha\sb s$ is non-zero,
we have $G\sb 1-G\sb 2\notin \IMb$,
which contradicts to~\eqref{eq:implication}.
Hence there
must exist
a point $(\tilde g, \tilde h)\in U$
in an arbitrary small neighborhood of $(g, h)$
such that~\eqref{eq:rhosurj}
holds.
Therefore $\rho$ is dominant.
\par
%
%
Finally we calculate
$\dim\Ker (d\alpha)\sb Q$.
By Proposition~\ref{prop:dimKerdalpha},
we see that 
$\dim\Ker (d\alpha)\sb Q$ is equal to 
$$
\dim \Fba  - \dim\Ha-s+
\dim (\Maa / (J\sb g\Maa +\Io h (\Mbb)))\sb 0.
$$
The fourth term
is equal to
$n+r-s$ by
Proposition~\ref{prop:equivconds}.
Since $n+r-2s=-m+2l$, we complete the proof of the assertion (2).
\end{sproof}
Now we are ready to the proof of Main Theorem.
\begin{mainproof}
For a locally closed analytic subspace $A$ of $\Ha$,
we denote by
$$
\begin{array}{ccc}
\tlZZ\sb{A} &\maprightsp{\alpha\sb A} &\XX\sb A \\
\mapdownleft{\pi\sb A} & &\mapdownright{\phi\sb A}\\
F\sb{A} &\maprightsb{\rho\sb A} &A\\
\end{array}
$$
the pull-back of the right square of the diagram~\eqref{diag:univ}
by $A\hookrightarrow \Ha$.
\par
There exists
a Zariski open dense subset
$\UU$ of $\Ha$
such that
$$
\begin{array}{ccc}
\tlZZ\sb{\UU} &\maprightsp{\alpha\sb{\UU}} &\XX\sb{\UU} \\
\mapdownleft{\pi\sb{\UU}} & &\mapdownright{\phi\sb{\UU}}\\
F\sb{\UU} &\maprightsb{\rho\sb{\UU}} &\UU\\
\end{array}
$$
is locally trivial
over $\UU$ in the category of topological spaces
and continuous maps,  that
$\phi\sb{\UU}$ is smooth,
and that $\rho\sb{\UU}$ is smooth or $F\sb{\UU}$ is empty. 
It is enough to show $F\sb{\bb} (X\sb b)\ne \emptyset$ and 
$\Im \psi\sb{\bb} (X\sb b)\supseteq V\sb m (X\sb b, \Z)$
for at least one point $b$ of $\UU$,
where $X\sb b$ denotes
the complete intersection corresponding to a point $b$ of $\UU$.
\par
%
%
By Proposition~\ref{prop:dims2} (1), $F\sb{\bb} (X\sb b)$ is non-empty for any $b\in \UU$.
\par
%
%
By the assumption of Main Theorem,
the morphism $\alpha |\Xi\sb o : \Xi\sb o \to \Gamma\sb o$
is dominant,
and hence, by Proposition~\ref{prop:uniqueirred},
there exists a unique irreducible component
$\Xi\sb o\sprime$ of $\Xi\sb o$
that is mapped dominantly onto $\Gamma\sb o$ by $\alpha|\Xi\sb o$.
Let $Q:=(o, \rep g, \rep f)$ be a general point of $\Xi\sb o\sprime$.
Then 
$\alpha (Q) =(o, \rep f)$ is
a general point of $\Gamma\sb o$.
By Corollary~\ref{cor:generalinGammao},
the point $o$ is the  only  singular point of $X\sb{\rep f}$,
and it  is a hypersurface singularity with non-degenerate Hessian.
In particular,
the image of $(d\phi)\sb{\alpha (Q)} : T\sb{\alpha (Q)} \XX \to T\sb{\rep f} \Ha$ is 
of codimension $1$ in $T\sb{\rep f} \Ha$.
On the other hand,
by Proposition~\ref{prop:dims2} (3),
the image of $(d\rho)\sb{\pi (Q)} : T\sb{\pi (Q)} \Fba \to T\sb{\rep f} \Ha$ is 
also of codimension $1$ in $T\sb{\rep f} \Ha$.
Hence there exists a smooth curve $C$ in $\Ha$
passing through $\rep f$
that  satisfies
\begin{equation}\label{eq:intTC}
\Im (d\phi)\sb{\alpha (Q)} \cap T\sb{\rep f} C=0 ,
\quad
\Im (d\rho)\sb{\pi (Q)} \cap T\sb{\rep f} C=0,
\end{equation}
and $C\cap \UU\ne \emptyset$.
We choose a sufficiently small
open unit disk
$\Delta $ in $C$
with the center $\rep f$,
and consider the following diagrams:
\begin{equation}\label{diag:overDelta}
\begin{array}{ccc}
\tlZZ\sb C &\maprightsp{\alpha\sb C} &\XX\sb C \\
\mapdownleft{\pi\sb C} & &\mapdownright{\phi\sb C}\\
F\sb C &\maprightsb{\rho\sb C} &C\\
\end{array}
\qquad\textrm{and}\qquad
\begin{array}{ccc}
\tlZZ\sbD &\maprightsp{\alpha\sbD} &\XX\sbD \\
\mapdownleft{\pi\sbD} & &\mapdownright{\phi\sbD}\\
F\sbD &\maprightsb{\rho\sbD} &\Delta.\\
\end{array}
\end{equation}
We can assume that $\Delta\sp{\times}:=\Delta\setminus \{\rep f\}$ is contained in $ \UU$.
By Lemma~\ref{lem:cut},
the analytic space
$\XX\sbD$ is smooth of dimension $m+1$.
Moreover,
the holomorphic map $\phi\sbD : \XX\sbD \to \Delta$
has only one critical point,
which is  the point  $(o, \rep f)$ on the central fiber $X\sb{\rep f}$,
and at which the  Hessian of $\phi\sb{\Delta}$ is  non-degenerate.
We select a point $b$ of $\UU$ from $\Delta\sp{\times}$.
Then 
we have a vanishing cycle $[\Sigma\sb b ]\in H\sb m (X\sb b, \Z)$,
unique up to sign, 
associated to the non-degenerate critical point $(o, \rep f)$ of $\phi\sb{\Delta}$.
It is known that 
$V\sb m (X\sb b, \Z)$
is generated by $[\Sigma\sb b ]$ as 
a module over the group ring $\Z [\pi\sb 1 (\UU, b)]$.
(See~\cite{Lamotke81}.)
On the other hand,
the image of the cylinder homomorphism
$\psi\sb{\bb} (X\sb b)$
is $\pi\sb 1 (\UU, b)$-invariant.
Therefore it is enough to show 
that the image of
$\psi\sb{\bb} (X\sb b)$
contains $[\Sigma \sb b]$.
\par
%
%
We put $O:=\pi (Q)\in \Fba$.
By Lemma~\ref{lem:cut} and~\eqref{eq:intTC},
the scheme $F\sb {C}$ in the left diagram of~\eqref{diag:overDelta}
is smooth at $O$,
and 
\begin{equation}\label{eq:dimF}
\dim\sb{O} F\sb {C}
=\dim \Fba -\dim \Ha +1.
\end{equation}
From the construction of $\tlZZ\sb C$, we see that 
$\Ker (d\alpha)\sb Q$
is contained 
in the subspace $T\sb Q \tlZZ\sb C$ of $T\sb Q \tlZZ$,
and that $\Ker (d\alpha)\sb Q$ coincides with 
$\Ker (d\alpha \sb {C} )\sb Q$.
Hence, by Proposition~\ref{prop:dims2} (2),
we have
\begin{equation}\label{eq:dimKeralpha}
\dim \Ker (d\alpha \sb {C})\sb Q = \dim \Fba -\dim \Ha -m+2l.
\end{equation}
Since $\tlZZ\sb C$
is a closed analytic subspace of $F\sb C \times \XX\sb C$
with $\pi\sb C$ and $\alpha\sb C$ being projections,
we have
\begin{equation}\label{eq:intzero}
\Ker (d\pi \sb {C})\sb Q \cap \Ker (d\alpha \sb {C})\sb Q=0
\end{equation}
in $T\sb Q \tlZZ\sb C$.
In particular,
the linear map $(d\pi\sb C)\sb Q : T\sb Q \tlZZ\sb C \to T\sb O  F\sb C$
maps $\Ker (d\alpha \sb {C})\sb Q$
isomorphically to a linear subspace of $T\sb{O} F\sb C$.
By the dimension counting~\eqref{eq:dimF} and~\eqref{eq:dimKeralpha},
this subspace
$$
(d\pi\sb C)\sb Q (\Ker (d\alpha \sb {C})\sb Q)\;\subset\; T\sb{O} F\sb C
$$
is of codimension $m-2l+1$.
Hence there exists a
closed subvariety $F\sprime\sb C$
of $F\sb C$
with dimension $m-2l+1$
that passes through $O$,
is
smooth at $O$, 
and
satisfies
\begin{equation}\label{eq:intzerodown}
T\sb{O} F\sprime\sb C \cap (d\pi\sb C)\sb Q (\Ker (d\alpha \sb{C})\sb Q)=0.
\end{equation}
We put 
$$
F\sprime\sbD:=F\sprime\sb C\cap F\sbD, \quad
\tlZZ\sprime\sb C:=\pi\sb C\inv (F\sprime\sb C)
\quand
\tlZZ\sprime\sbD:=\pi\sbD\inv (F\sprime\sbD),
$$ 
and let
\begin{equation}\label{diag:overDeltasprime}
\begin{array}{ccc}
\tlZZ\sprime\sb C &\maprightsp{\alpha\sprime\sb C} &\XX\sb C \\
\mapdownleft{\pi\sprime\sb C} & &\mapdownright{\phi\sb C}\\
F\sprime\sb C &\maprightsb{\rho\sprime\sb C} &C\\
\end{array}
\qquad\textrm{and}\qquad
\begin{array}{ccc}
\tlZZ\sprime\sbD &\maprightsp{\alpha\sprime\sbD} &\XX\sbD \\
\mapdownleft{\pi\sprime\sbD} & &\mapdownright{\phi\sbD}\\
F\sprime\sbD &\maprightsb{\rho\sprime\sbD} &\Delta \\
\end{array}
\end{equation}
be the restriction of the diagrams~\eqref{diag:overDelta}.
The right diagram of~\eqref{diag:overDeltasprime}
is the pull-back of the left diagram of~\eqref{diag:overDeltasprime}
by $\Delta\hookrightarrow C$.
\par
%
%
Since the fiber of $\pi$ passing through $Q$ is smooth at $Q$
by the definition of $\Xi\sb o$,
the holomorphic map $\pi\sprime\sbD$ is also smooth at $Q$.
Moreover,
from~\eqref{eq:intzero} and~\eqref{eq:intzerodown},
we have
\begin{equation*}
\begin{split}
\Ker (d\alpha\sprime\sbD)\sb Q\;\;=\;\;&T\sb Q \tlZZ\sprime\sbD \cap \Ker (d\alpha\sbD)\sb Q\\ 
\;\;=\;\;&(d\pi\sb C)\sb Q\inv (T\sb O F\sb C\sprime) \cap \Ker (d\alpha\sb C)\sb Q
\;\;=\;\;0.
\end{split}
\end{equation*}
Therefore $\alpha\sprime\sbD$ is
an immersion at $Q$.
We have $\dim F\sprime\sb{\Delta} =m-2l+1$.
Note that  $ H\sb m (X\sb b , \Z)$
is torsion free.
Hence
the right diagram of~\eqref{diag:overDeltasprime}
satisfies all the conditions required in Theorem~\ref{thm:vccyl} (2).
We put
$$
F\sbD \sprime (X\sb b) :=\rho\sp{\prime\;-1}\sbD (b), \quad
Z\sbD \sprime (X\sb b):=\pi\sp{\prime\;-1	}\sbD(F\sprime\sbD (X\sb b)),
$$
and consider the family
\begin{equation}\label{eq:sprimeb}
\qquad\begin{array}{ccc}
 Z\sbD \sprime (X\sb b)&\maprightsp{} &X\sb b \\
\mapdownleft{} && \\
F\sbD \sprime (X\sb b) &&
\end{array}
\end{equation}
of $l$-dimensional closed analytic subspaces  of $X\sb b$.
By Theorem~\ref{thm:vccyl}, the image of the cylinder homomorphism
$$
\map{\psi\sprime\sb{\bb} (X\sb b)}{H\sb{m-2l} (F\sbD \sprime (X\sb b), \Z)}{H\sb m (X\sb b, \Z)}
$$
associated with  the family~\eqref{eq:sprimeb}
contains the vanishing cycle $[\Sigma\sb b] \in H\sb m (X\sb b, \Z)$.
By the construction,
$\psi\sprime\sb{\bb} (X\sb b)$
is the composite
of the homomorphism
$$
H\sb{m-2l} (F\sbD\sprime (X\sb b) , \Z) \to H\sb{m-2l}(F\sb{\bb} (X\sb b), \Z)
$$
induced from the inclusion
$F\sbD\sprime (X\sb b) \hookrightarrow F\sb{\bb} (X\sb b)$
and the original cylinder homomorphism  $\psi\sb{\bb} (X\sb b)$.
Hence
the image of $\psi\sb{\bb}(X\sb b)$ contains
$[\Sigma\sb b]$.
\end{mainproof}
We put $F\sb C\sprime (X\sb b):=\rho\sb C\sp{\prime\;-1} (b)$,
and 
let $F\sb C\spprime (X\sb b)$ be the union of irreducible
components of  $F\sb C\sprime (X\sb b)$
with dimension $m-2l$.
Then $F\sb C\spprime (X\sb b)$
contains 
an $(m-2l)$-dimensional sphere 
representing the vanishing cycle $[\sigma\sb b]\in H\sb{m-2l} (F\sbD\sprime (X\sb b) , \Z)$ 
associated to the non-degenerate critical point $O$ of $\rho\sprime\sbD$.
Let $T$ be the Zariski closure of 
$\alpha\sprime\sb C(\pi\sp{\prime\;-1}\sb C(F\sb C\spprime (X\sb b)))$
in $X\sb b$.
Then $T$ is of dimension  $m-l$, and 
$[\Sigma\sb b]=\pm \psi\sprime\sb{\bb} (X\sb b) ([\sigma\sb b])$
is represented by a topological cycle whose support is contained in $T$.
Therefore we obtain the following:
\begin{scorollary}\label{cor:toGHC}
Suppose that $(n,\aa, \bb)$ satisfies the conditions of Main Theorem.
Let $X$ be a general complete intersection of multi-degree $\aa$ in $\P\sp n$.
Then 
every  vanishing cycle
of $X$
is represented by a topological cycle
whose support is contained in
a Zariski closed subset 
of  $X$ with codimension $l$.
\qed
\end{scorollary}
\section{Gr\"obner bases method}\label{sec:GBmethod}
Suppose  we are given a triple $(n, \aa,\bb)$
that satisfies 
the  conditions~\eqref{eq:num1} and~\eqref{eq:num2}
of Main Theorem.
We will describe a method to determine
whether this triple satisfies the second condition of Main Theorem.
\par
First we choose a prime integer $p$,
and put
$$
R\spp:=\F\sb p [x\sb 0, \dots, x\sb n].
$$
We define graded $R\spp$-modules
$\Maa\spp$, $\Mbb\spp$, $\Naa\spp$,
and ideals $I\sbo\spp$, $J\sb g\spp$ of $R\spp$ in the same way
as in  \S\ref{sec:univ} except for the coefficient field.
We generate
an element 
$g=\transpose{ (g\sb 1, \dots,g\sb s)}$ of $  (\Io\spp \Mbb\spp)\sb 0$
and
a homomorphism 
$h=(h\sb{ij})\in \Hom (\Mbb\spp,\Naa\spp)\sb 0$
in a random way.
Then we can calculate
\begin{equation}\label{eq:dimp}
\dim\sb{\F\sb p} (\Maa\spp/(J\sb g\spp \Maa\spp +\Io\spp h (\Mbb\spp)) )\sb 0
\end{equation}
by means of Gr\"obner bases.
If this dimension is  $\le n+r-s$,
then the condition (iv) of
Proposition~\ref{prop:equivconds}
is fulfilled,
because this condition is an open condition.
Hence the morphism $\alpha |\Xi\sb o: \Xi\sb o \to \Gamma\sb o$ is dominant.
%
%
%
%
%
%
\section{Application of   a theorem
of Debarre and Manivel}\label{sec:DMmethod}
From now on,
we use the following terminology.
A {\em sequence}
always means  a finite non-decreasing sequence
of positive integers.
For a sequence $\aa$,
let $\min (\aa)$ and $\max (\aa)$
be the first and the last elements of $\aa$, respectively, 
and let $|\aa|$ denote the length of $\aa$.
Let $\aa\sprime$ be another sequence.
We denote by
$\aa\uplus\aa\sprime$ the sequence of
length $|\aa|+|\aa\sprime|$
obtained 
by re-arranging the conjunction $(\aa, \aa\sprime)$
into the non-decreasing order.
For an integer $a\ge 2$,
we define 
$(a)!$ to be the sequence $(2, \dots, a)$
of length $a-1$,
and
for a sequence $\aa=(a\sb 1, \dots, a\sb r)$
with $\min (\aa)\ge 2$,
we put
$$
\aa !:=(a\sb 1)!\uplus \dots \uplus (a\sb r)!\;.
$$
We sometimes write a sequence by indicating 
the number of repetition of each integer
in the sequence by a superscript.
For example, we have
$(2,3,3,4)!=(2,2,2,2,3,3,3,4)=(2\sp 4, 3\sp 3, 4)$.
\par
%
%
Let $n$ and $\ell$ be positive integers,
and $\aa=(a\sb 1, \dots, a\sb r)$
a sequence.
According to~\cite{DebarreManivel98},
we put
$$
\delta (n, \aa,\ell):=(\ell+1)(n-\ell) -\sum\sb{i=1}\sp r \binomial{a\sb i+\ell}{\ell},
$$
and $\delta\sb{-} (n,\aa, \ell):=\min\{\delta (n,\aa, \ell), n-2\ell-|\aa|\}$.
\begin{stheorem}[\cite{DebarreManivel98}, Th\'eor\`eme 2.1]\label{thm:DM}
A general complete intersection of multi-degree $\aa$ in $\P\sp n$
contains an $\ell$-dimensional linear subspace
if and only if $\delta\sb{-} (n, \aa, \ell)\ge 0$.
\noproof
\end{stheorem}
\begin{stheorem}\label{thm:numlin}
Let $\aa=(a\sb 1, \dots, a\sb r)$ be a sequence
satisfying $\min (\aa) \ge 2$ and $\sum\sb{i=1}\sp r  a\sb i \le n$.
Let $\aa\sprime$ be a sub-sequence
of $\aa$ such that $\max (\aa\sprime)=\max (\aa)$,
and let $\aa\spprime$ be the complement to $\aa\sprime$ in $\aa$.
{\rm (}%
When $\aa\sprime=\aa$, 
$\aa\spprime$  is the empty sequence.%
{\rm )}
Suppose that a positive integer $\lambda$
satisfies the following\itcolon
\begin{equation}\label{eq:numcondlin}
\delta\sb{-} (n-|\aa\sprime|, \aa\sprime !\,, \lambda-1)\ge 0,
\quad
|\aa\spprime|<\lambda 
\quand
n-r>2 (\lambda -|\aa\spprime |).
\end{equation}
We put 
$\bb:=(1\sp{n-\lambda}) \uplus \aa\spprime$.
Then 
$F\sb{\bb} (X)$ is non-empty
for a general complete intersection $X$ of multi-degree $\aa$ in $\P\sp n$,
and 
the image of the cylinder homomorphism $\psi\sb{\bb} (X)$
contains 
$V\sb m (X, \Z)$.
\end{stheorem}
\begin{sproof}
Note that $l=n-|\bb|$ is equal to $\lambda-|\aa\spprime|$.
Since $\max (\aa\sprime)=\max (\aa)$,
we can assume that $a\sb r$ is a member of $\aa\sprime$.
Let $f=\transpose{(f\sb 1, \dots,f\sb r)}$ be a general element  of $\IN$,
and let
$Y\sb i$ be the hypersurface of degree $a\sb i$
defined by $f\sb i=0$.
We put 
%
%
$$
X\sprime :=\bigcap \sb{a\sb i \in \aa\sprime} Y\sb i\quand
X\spprime :=\bigcap \sb{a\sb i \in \aa\spprime} Y\sb i.
$$
By Proposition~\ref{prop:Gammao}, $X\sprime$ is a general member
of the family of   complete intersections of
multi-degree $\aa\sprime$
possessing a singular point at $o$.
By means of the projection with the center $o$,
we see that $X\sprime$ contains 
a linear subspace of dimension $\ell>0$ that passes through $o$
if and only if
a general complete intersection of multi-degree $\aa\sprime !$
in $\P\sp{n-|\aa\sprime|}$
contains an $(\ell-1)$-dimensional linear subspace.
By Theorem~\ref{thm:DM},
the first condition 
of~\eqref{eq:numcondlin}
implies that $X\sprime$ contains a linear subspace $\Lambda$
of dimension $\lambda$ passing through $o$.
In particular,
we have
$\lambda < n-|\aa\sprime|$.
Using this inequality
and the second and the third conditions
of~\eqref{eq:numcondlin},
we can easily check that
$(n, \aa,\bb)$
satisfies
the  conditions~\eqref{eq:num1} and~\eqref{eq:num2}
in Main Theorem.
\par
We put $Z:=\Lambda\cap X\spprime$.
Then $Z$ is a complete intersection of multi-degree $\bb$
contained in $X\sb{\rep f} =X\sprime\cap X\spprime$
and passing through $o$.
Moreover,
since the polynomials $f\sb i$ $(a\sb i \in \aa\spprime)$
are general with respect to
$\Lambda$, $Z$ is smooth.
Thus the second condition of Main Theorem is also
satisfied.
\end{sproof}
\section{The generalized Hodge conjecture for complete intersections}\label{sec:GHC}
\begin{table}[bt]
{\small
\newcommand{\stample}[2]{
\noline
\strut  &\;${#1}$\; &&\;${#2}$\;\hfill& \cr
\noline
}
\newcommand{\twoline}{%
\multispan{5}\hrulefill \cr
}
\newcommand{\noline}{%
height 1pt &&&&\cr
}
\vbox{\offinterlineskip
\halign{ 
\vrule # & # & \vrule # &\hskip 3pt   #&\vrule #  \cr 
\twoline
\noline
\strut & \; $n$\; & & \hfill $\aa$\phantom{$[]$} \hfill & \cr
\noline
\twoline
\noline
\twoline
\stample{6}{
({3})}
\twoline
\stample{7}{
({3})}
\twoline
\stample{8}{
({2}, {3}), ({3}), ({4})}
\twoline
\stample{9}{
({2}, {3}), ({3}), ({3}\sp{2}), ({4})}
\twoline
\stample{10}{
({2}, {3}), ({3}), ({3}\sp{2}), ({4})}
\twoline
\stample{11}{
({2}\sp{2}, {3}), ({2}, {4}), ({3}), ({3}\sp{2}), ({4}), ({5})}
\twoline
\stample{12}{
({2}\sp{2}, {3}), ({2}, {3}), ({2}, {3}\sp{2}), ({2}, {4}), ({3}), ({3}\sp{2}), ({3}\sp{3}), ({5})}
\twoline
\stample{13}{
({2}\sp{3}, {3}), ({2}\sp{2}, {4}), ({2}, {3}), ({2}, {3}\sp{2}), ({2}, {4}), ({3}), ({3}\sp{2}), ({3}\sp{3}), ({3}, {4}), ({5})}
\twoline
\stample{14}{
({2}\sp{3}, {3}), ({2}\sp{2}, {4}), ({2}, {5}), ({3}), ({3}\sp{2}), ({3}\sp{3}), ({3}, {4}), ({4}), ({4}\sp{2}), ({5})}
\twoline
\stample{15}{
({2}\sp{2}, {3}), ({2}\sp{2}, {3}\sp{2}), ({2}\sp{2}, {4}), ({2}, {3}), ({2}, {3}\sp{2}), ({2}, {3}\sp{3}), ({2}, {3}, {4}), ({2}, {5}), ({3}), ({3}\sp{2}),
}
\stample{}{({3}\sp{3}), ({3}\sp{4}), ({4}), ({4}\sp{2}), ({6})}
\twoline
\stample{16}{
({2}\sp{4}, {3}), ({2}\sp{3}, {4}), ({2}\sp{2}, {5}), ({2}, {3}), ({2}, {3}\sp{2}),({2}, {3}\sp{3}), ({2}, {3}, {4}), ({2}, {5}), ({3}), ({3}\sp{2}),}
\stample{}{ 
({3}\sp{3}), ({3}\sp{4}), ({3}, {5}), ({6})}
\twoline
\stample{17}{
({2}\sp{3}, {4}), ({2}\sp{2}, {5}), ({2}, {4}), ({2}, {4}\sp{2}), ({2}, {6}), ({3}), ({3}\sp{2}), ({3}\sp{3}), ({3}\sp{4}), ({3}\sp{2}, {4}), ({3}, {5}), ({6})}
\twoline
\stample{18}{
({2}\sp{5}, {3}), ({2}\sp{4}, {4}), ({2}\sp{2}, {3}, {4}), ({2}\sp{2}, {5}), ({2}, {3}, {5}), ({2}, {6}), ({4}, {5})}
\twoline
\stample{19}{
({2}\sp{4}, {4}), ({2}\sp{3}, {5}), ({2}, {3}), ({2}, {3}\sp{2}), ({2}, {3}\sp{3}), ({2}, {3}\sp{4}), ({2}, {3}, {5}), ({2}, {6}), ({3}), ({3}\sp{2}),}
\stample{}{
({3}\sp{3}), ({3}\sp{4}), ({3}\sp{5}), ({3}, {6}), ({7})}
\twoline
\stample{20}{
({2}\sp{3}, {3}, {4}), ({2}\sp{3}, {5}), ({2}\sp{2}, {6}), ({3}), ({3}\sp{2}), ({3}\sp{3}), ({3}\sp{4}), ({3}\sp{5}), ({3}\sp{2}, {5}), ({3}, {6}), ({7})}
\twoline
\stample{21}{
({2}\sp{5}, {4}), ({2}\sp{4}, {5}), ({2}\sp{2}, {3}, {5}), ({2}\sp{2}, {6}), ({2}, {7})}
\twoline
\stample{22}{
({2}\sp{4}, {5}), ({2}\sp{3}, {6}), ({2}, {3}, {6}), ({2}, {7})}
\twoline
\stample{23}{
({2}\sp{6}, {4}), ({2}\sp{3}, {3}, {5}), ({2}\sp{3}, {6}), ({3}), ({3}\sp{2}), ({3}\sp{3}), ({3}\sp{4}), ({3}\sp{5}), ({3}\sp{6}), ({3}, {7}), ({8})}
\twoline
\stample{24}{
({2}\sp{5}, {5}), ({2}\sp{2}, {3}, {6}), ({2}\sp{2}, {7})}
\twoline
\stample{25}{
({2}\sp{4}, {6})}
\twoline
\stample{26}{
({2}\sp{6}, {5}), ({2}\sp{3}, {7})}
\twoline
\stample{27}{
({2}\sp{5}, {6})}
\twoline
}
}
\vskip 5pt
}
\caption{The $148$ pairs}
\label{table:A}
\end{table}
\begin{table}[bt]
{\small
\newcommand{\stample}[3]{\strut  &\,${#1}$\, &&\;${#2}$\;&&\;${#3}$\;& \cr}
\newcommand{\sstample}[2]{\strut  & &&\;${#1}$\;&&\;${#2}$\;& \cr}
\newcommand{\twoline}{%
\noline
&& \multispan{5}\hrulefill \cr
\noline
}
\newcommand{\threeline}{%
\noline
\multispan{7}\hrulefill \cr
\noline
}
\newcommand{\upthinthreeline}{%
\multispan{7}\hrulefill \cr
\noline
}
\newcommand{\downthinthreeline}{%
\noline
\multispan{7}\hrulefill \cr
}
\newcommand{\noline}{%
height 1pt &&&&&&\cr
}
\setbox1=\hbox{\vbox{\offinterlineskip
\halign{ 
\vrule # & # & \vrule # &  #&\vrule # & # &\vrule # \cr 
\upthinthreeline
\strut &\;${n}$\; &&\hfill\;${\aa}$\;\hfill&&\hfill\;${\bb}$\;\hfill& \cr
\downthinthreeline
\threeline
  \stample{10}{(2 \sp {2}, 3)}{(1 \sp {7},2)}
\threeline
  \stample{11}{(2, 3)}{(1 \sp {7},2)}
\twoline
 \sstample{(2, 3 \sp {2})}{(1 \sp {7},2,3)}
\threeline
  \stample{12}{(2 \sp {3}, 3)}{(1 \sp {7},2 \sp {3})}
\threeline
  \stample{13}{(2 \sp {2}, 3)}{(1 \sp {8},2 \sp {2})}
\twoline
 \sstample{(2 \sp {2}, 3 \sp {2})}{(1 \sp {8},2 \sp {2},3)}
\threeline
  \stample{14}{(2 \sp {2}, 3)}{(1 \sp {10},2)}
\twoline
 \sstample{(2 \sp {2}, 3 \sp {2})}{(1 \sp {9},2 \sp {2},3)}
\twoline
 \sstample{(2, 3)}{(1 \sp {9},2)}
\twoline
 \sstample{(2, 3 \sp {2})}{(1 \sp {9},2,3)}
\twoline
 \sstample{(2, 3 \sp {3})}{(1 \sp {9},2,3 \sp {2})}
\downthinthreeline
}}}
\setbox2=
\hbox{
\vbox{\offinterlineskip
\halign{ 
\vrule # & # & \vrule # &  #&\vrule # & # &\vrule # \cr 
\upthinthreeline
\strut &\;${n}$\; &&\hfill\;${\aa}$\;\hfill&&\hfill\;${\bb}$\;\hfill& \cr
\downthinthreeline
\threeline
  \stample{15}{(2 \sp {4}, 3)}{(1 \sp {9},2 \sp {4})}
\threeline
  \stample{16}{(2 \sp {3}, 3)}{(1 \sp {10},2 \sp {3})}
\twoline
 \sstample{(2 \sp {3}, 3 \sp {2})}{(1 \sp {10},2 \sp {3},3)}
\threeline
  \stample{17}{(2 \sp {3}, 3)}{(1 \sp {13},2)}
\twoline
 \sstample{(2 \sp {3}, 3 \sp {2})}{(1 \sp {11},2 \sp {3},3)}
\twoline
 \sstample{(2 \sp {2}, 3)}{(1 \sp {11},2 \sp {2})}
\twoline
 \sstample{(2 \sp {2}, 3 \sp {2})}{(1 \sp {11},2 \sp {2},3)}
\twoline
 \sstample{(2 \sp {2}, 3 \sp {3})}{(1 \sp {11},2 \sp {2},3 \sp {2})}
\threeline
  \stample{18}{(2 \sp {2}, 3)}{(1 \sp {13},2)}
\twoline
 \sstample{(2 \sp {2}, 3 \sp {2})}{(1 \sp {13},2,3)}
\twoline
 \sstample{(2 \sp {2}, 3 \sp {3})}{(1 \sp {13},2,3 \sp {2})}
\downthinthreeline
}}}
\setbox3=
\hbox{
\vbox{\offinterlineskip
\halign{ 
\vrule # & # & \vrule # &  #&\vrule # & # &\vrule # \cr 
\upthinthreeline
\strut &\;${n}$\; &&\hfill\;${\aa}$\;\hfill&&\hfill\;${\bb}$\;\hfill& \cr
\downthinthreeline
\threeline
  \stample{18}{(2, 3)}{(1 \sp {12},2)}
\twoline
 \sstample{(2, 3 \sp {2})}{(1 \sp {12},2,3)}
\twoline
 \sstample{(2, 3 \sp {3})}{(1 \sp {12},2,3 \sp {2})}
\twoline
 \sstample{(2, 3 \sp {4})}{(1 \sp {12},2,3 \sp {3})}
\threeline
 \stample{19}{(2 \sp {4}, 3)}{(1 \sp {12},2 \sp {4})}
\twoline
 \sstample{(2 \sp {4}, 3 \sp {2})}{(1 \sp {12},2 \sp {4},3)}
\threeline
  \stample{20}{(2 \sp {3}, 3)}{(1 \sp {13},2 \sp {3})}
\twoline
 \sstample{(2 \sp {3}, 3 \sp {2})}{(1 \sp {13},2 \sp {3},3)}
\twoline
 \sstample{(2 \sp {3}, 3 \sp {3})}{(1 \sp {13},2 \sp {3},3 \sp {2})}
\downthinthreeline
}
}
}
\vbox to \ht1  
{\hbox
  {
  \raise \ht2 \hbox{\raise \ht3 \copy1}
  \hskip 2pt 
  \raise \ht1 \hbox{\raise \ht3 \copy2}
  \hskip 2pt 
  \raise \ht1 \hbox{\raise \ht2 \copy3}
  }
}
\vskip 5pt
}
\caption{Examples of triples obtained by Gr\"obner bases method}
\label{table:B}
\end{table}
Suppose we are given  a pair $(n, \aa)$
satisfying  $\min (\aa) \ge 2$ and $\sum a\sb i \le n$.
We put
$k:=[(n-\sum a\sb i)/\max (\aa)]+1$.
The
Hodge structure
of the middle cohomology group
$H\sp m (X)$
of a general complete intersection $X$
of multi-degree $\aa$
in $\P\sp n$ satisfies~\eqref{eq:HS}.
We will investigate the consequence of the generalized Hodge conjecture
that there should exist a Zariski closed subset $T$ of $X$
with codimension $k$ such that
every element of $H\sb m (X, \Q)$
is represented by a topological cycle 
whose support is contained in $T$.
Note that $H\sb m (X, \Q)$ is generated by vanishing cycles
and, if $m$ is even,
the homology class of an  intersection of $X$ and a linear subspace of $\P\sp n$.
Hence, by Corollary~\ref{cor:toGHC},
this consequence is verified if we can find 
$\bb$ with the following properties:
\begin{equation}\label{eq:goal}
\begin{split}
&l:=n-|\bb|=k\quand\\
&\textrm{$(n, \aa, \bb)$ satisfies the assumptions of Main Theorem.}
\end{split}
\end{equation}
In the following, we assume $m>2k$.
This inequality  $m>2k$ fails to hold if and only if 
 $m\le 2$ or $\aa=(2)$ or
($\aa=(2, 2)$ and $m$ even).
In these cases, the Hodge conjecture
has been already proved.
\begin{sproposition}[\cite{Shimada90a}, \cite{Shimada91}]\label{prop:koneandquads}
\par
{\rm (1)}
If $k=1$,
then $\bb=(1\sp{n-1})$
satisfies~\eqref{eq:goal}.
\par
{\rm (2)}
If $\aa=(2\sp r)$,
then
$\bb=(1\sp{n-[n/2]}, 2\sp{r-1})$
satisfies~\eqref{eq:goal}.
\end{sproposition}
\begin{sproof}
Put $\aa\sprime=\aa$ in the case (1) and  $\aa\sprime=(2)$ in the  case  (2), and 
apply Theorem~\ref{thm:numlin}.
\end{sproof}
In these cases,
the consequence of the generalized Hodge conjecture is verified in any dimension.
\par
We have made an exhaustive
search in $n\le 40$,
and found $148$ pairs $(n, \aa)$
that are not covered by 
Proposition~\ref{prop:koneandquads},
but for which Theorem~\ref{thm:numlin}
yields $\bb$ satisfying~\eqref{eq:goal}
by taking an appropriate sub-sequence $\aa\sprime$.
We list up these $(n, \aa)$ in Table~\ref{table:A}.
No such $(n, \aa)$ are found in $n>27$.
Even if $(n, \aa)$ does not appear
in Table~\ref{table:A},
the calculation of the dimension~\eqref{eq:dimp} by 
Gr\"obner bases
sometimes
gives us $\bb$ with~\eqref{eq:goal}.
Examples of these  $(n, \aa, \bb)$
in $n\le 20$ are given in Table~\ref{table:B}.
From these results,
we can find $\bb$ with~\eqref{eq:goal}
for any $(n, \aa)$ with $n\le 9$.
When $n=10$,
$\aa=(2,4)$ and $\aa=(5)$
appear in neither Tables~\ref{table:A} nor~\ref{table:B}.
\par
As a closing remark,
let us return to the classical example
of cubic threefolds (\cite{ClemensGriffiths72}).
Our method shows that,
not only the family of lines $\bb=(1\sp 3)$,
but  also the family of curves
with $\bb=(1\sp 2, 2)$ or $(1, 2\sp 2)$ or $(2\sp 3)$
give a surjective cylinder homomorphism
on the middle homology group $H\sb 3 (X, \Z)$
of a general cubic threefold $X$.
%
%
%
%
%
\par
\bigskip
\noindent
{\small
{\bf Acknowledgment.}
Part of this work was done
during the author's stay at Nagoya University in December 2000.
He would like to thank
Professors
Shigeyuki Kondo and 
Hiroshi Saito
for their warm hospitality.
}

\end{document}